\documentclass{amsart}
\usepackage{amssymb}
\usepackage[all]{xy}

\theoremstyle{definition}
\newtheorem{defn}{Definition}[section]
\theoremstyle{remark}
\newtheorem{rem}[defn]{Remark}
\newtheorem{rems}[defn]{Remarks}
\newtheorem{example}[defn]{Example}
\newtheorem{examples}[defn]{Examples}
\theoremstyle{plain}
\newtheorem{prop}[defn]{Proposition}
\newtheorem{cor}[defn]{Corollary}
\newtheorem{thm}[defn]{Theorem}
\newtheorem{lemma}[defn]{Lemma}

\newcommand{\rank}{\mathrm{rk}}
\newcommand{\Cohstack}{\mathcal{C}oh}
\newcommand{\Bunstack}{\mathcal{B}un}
\newcommand{\Vect}{\mathcal{V}ect}

\newcommand{\schemes}{\mathsf{Schemes}}

\newcommand{\sets}{\mathsf{Sets}}
\newcommand{\Aa}{\mathbb{A}}
\newcommand{\Parstack}{\mathcal{P}ar}

\newcommand{\Extstack}{\mathcal{E}xt}
\newcommand{\Gr}{\mathrm{Gr}}
\newcommand{\GL}{\mathrm{GL}}
\newcommand{\Br}{\mathrm{Br}}

\newcommand{\cohom}{\mathrm{H}}
\newcommand{\Hom}{\mathrm{Hom}}

\newcommand{\Ext}{\mathrm{Ext}}
\newcommand{\End}{\mathrm{End}}
\newcommand{\Hombdl}{\mathit{Hom}}
\newcommand{\Endbdl}{\mathit{End}}
\newcommand{\Aut}{\mathrm{Aut}}

\newcommand{\longto}[1][]{\stackrel{#1}{\longrightarrow}}
\renewcommand{\to}[1][]{\stackrel{#1}{\rightarrow}}
\newcommand{\Quot}{\mathrm{Quot}}
\newcommand{\Tot}{\mathrm{Tot}}
\newcommand{\Spec}{\mathrm{Spec}}

\newcommand{\Pic}{\mathrm{Pic}}
\newcommand{\stackM}{\mathcal{M}}
\newcommand{\Buncoarse}{\mathfrak{Bun}}
\newcommand{\Parcoarse}{\mathfrak{Par}}
\newcommand{\Levcoarse}{\mathfrak{Lev}}
\newcommand{\Cohcoarse}{\mathfrak{Coh}}
\newcommand{\BNcoarse}{\mathfrak{BN}}
\newcommand{\Deccoarse}{\mathfrak{Dec}}
\newcommand{\coarseU}{\mathfrak{U}}

\newcommand{\coarseM}{\mathfrak{M}}
\newcommand{\coarseL}{\mathfrak{L}}
\renewcommand{\O}{\mathcal{O}}
\newcommand{\univ}{\mathrm{univ}}
\newcommand{\stab}{\mathrm{stab}}

\newcommand{\Suniv}{\mathcal{S}^\univ}
\newcommand{\Euniv}{\mathcal{E}^\univ}

\newcommand{\stackE}{\mathcal{E}}
\newcommand{\stackF}{\mathcal{F}}
\newcommand{\stackL}{\mathcal{L}}
\newcommand{\stackU}{\mathcal{U}}
\newcommand{\stackV}{\mathcal{V}}
\newcommand{\stackW}{\mathcal{W}}
\newcommand{\stackS}{\mathcal{S}}

\newcommand{\stackD}{\mathcal{D}}
\newcommand{\integers}{\mathbb{Z}}
\newcommand{\naturals}{\mathbb{N}}
\newcommand{\complexnums}{\mathbb{C}}
\newcommand{\reals}{\mathbb{R}}
\newcommand{\rationals}{\mathbb{Q}}
\newcommand{\dual}{\mathrm{dual}}
\newcommand{\coker}{\mathrm{coker}}
\newcommand{\im}{\mathrm{im}}
\newcommand{\id}{\mathrm{id}}
\newcommand{\pr}{\mathrm{pr}}
\newcommand{\op}{\mathrm{op}}
\newcommand{\hcf}{\mathrm{hcf}}
\newcommand{\Gm}{\mathbb{G}_m}
\newcommand{\BGm}{\mathrm{B}\mathbb{G}_m}
\newcommand{\BGL}{\mathrm{BGL}}

\begin{document}
\bibliographystyle{plain}

\title[Rationality and Poincar{\'e} families]{Rationality and Poincar\'{e} families\\for vector bundles with
extra structure\\on a curve}
\author[N. Hoffmann]{Norbert Hoffmann}
\address{Mathematisches Institut der Universit\"at\\Bunsenstr. 3--5\\37073 G\"ottingen\\Germany}
\curraddr{School of Mathematics\\Tata Institute of Fundamental Research\\Homi Bhabha Road\\Mumbai 400005\\
India}
\email{hoffmann@uni-math.gwdg.de}

\begin{abstract}
  Iterated Grassmannian bundles over moduli stacks of vector bundles on a curve are shown to be
  birational to an affine space times a moduli stack of degree $0$ vector bundles, following the method of
  King and Schofield. Applications include the birational type of some Brill-Noether loci, of moduli schemes
  for vector bundles with parabolic structure or with level structure and for A. Schmitt's decorated vector
  bundles. A further consequence concerns the existence of Poincar\'{e} families on finite coverings of the
  moduli schemes.
\end{abstract}

\maketitle

\section*{Introduction}
Let $C$ be a smooth projective algebraic curve of genus $g \geq 2$, say over an algebraically closed field
$k$, and let $L$ be a line bundle on $C$. A. King and A. Schofield \cite{king-schofield} have proved that the
coarse moduli scheme $\Buncoarse_{r, L}$ of stable vector bundles $E$ on $C$ with $\rank( E) = r$ and $\det(
E) \cong L$ is rational if the highest common factor $h$ of $r$ and $\deg(L)$ is $1$. The present paper
generalises this result to vector bundles with some extra structure, e.\,g.\ with parabolic structure in the
sense of \cite{mehta-seshadri}. In that case, we obtain that the moduli scheme is rational if the highest
common factor $h$ of rank, degree and all multiplicities is $1$, improving results of \cite{boden-yokogawa}.

This highest common factor $h$ also governs the existence of Poincar\'{e} families: Standard methods
\cite[Chapter 4, \S5]{newstead} construct e.\,g.\ a Poincar\'{e} family on the coarse moduli scheme of stable
parabolic vector bundles if $h = 1$, whereas Ramanan \cite{ramanan} has proved that there is no Poincar\'e
family on any open subscheme of $\Buncoarse_{r, L}$ if $h \neq 1$. We provide here a
common source for such results on the rationality of coarse moduli schemes \emph{and} on the existence of
Poincar\'{e} families on them: Our main theorem \ref{mainthm} states that the moduli \emph{stack} of fixed
determinant vector bundles with extra structure is birational to an affine space $\Aa^s$ times a moduli stack
$\Bunstack_{h, L_0}$ of rank $h$ vector bundles with fixed determinant $L_0$ of degree $0$.

The proof of this theorem requires to carefully keep track of something as inconspicuous as the scalar
automorphisms of the vector bundles $E$ on $C$. In some sense, this is already implicit in the King-Schofield
proof, namely in their notion of weight \cite[p. 526]{king-schofield} and in their Brauer class $\psi_{r, d}$
-- which can be interpreted as the obstruction against a Poincar\'{e} family. But the (admittedly more
abstract) stack language helps to clarify and generalise here, for example to arbitrary infinite base field
$k$ (as long as $C$ has a $k$-rational point). Apart from that, we basically follow the method of
\cite{king-schofield}: The moduli stack $\Bunstack_{r, L}$ is shown to be birational to a Grassmannian bundle
over $\Bunstack_{r_1, L_1}$ for some $r_1 < r$ if $r$ does not divide $\deg( L)$; then induction is used.

A new ingredient in our proof is the refinement \ref{ramanan} of Ramanan's theorem mentioned above; it allows
to extend our results to vector bundles with extra structures parameterised by iterated Grassmannian bundles
over $\Bunstack_{r, L}$, and it also yields that finite coverings of the coarse moduli scheme can only admit
Poincar\'{e} families if their degree is a multiple of $h$, cf. corollary \ref{no_poincare}.

The structure of this text is as follows: Section \ref{stacks} presents the moduli stacks that this paper
deals with, and collects some basic information about them. In section \ref{gerbes}, we introduce the notions
of $\Gm$-stack and of $\Gm$-gerbe in order to systematically distinguish `scalar automorphisms'. Section
\ref{vectorbundles} and section \ref{Grassmannians} deal with vector bundles and the associated Grassmannian
bundles on such stacks, respectively. Finally, section \ref{mainresults} contains statement and proof of the
main theorem as well as some consequences and examples, including vector bundles with parabolic structure,
Brill-Noether loci and A. Schmitt's decorated vector bundles.

I wish to thank J. Heinloth for numerous explanations and fruitful discussions about these moduli stacks, in
particular for his help in proving \ref{ramanan} and \ref{gentriv}. I also thank the Tata Institute of
Fundamental Research in Bombay for support and hospitality while this text was written.

\section{Moduli stacks of vector bundles on a curve} \label{stacks}
Let $k$ be an infinite field, and let $C$ be a geometrically irreducible, smooth projective curve over $k$ of
genus $g \geq 2$ which has a rational point $P \in C( k)$. This paper deals with moduli spaces of vector
bundles $E$ on $C$ from a birational point of view.
\begin{rem}
  Usually, such questions have been studied over algebraically closed fields $k$ or even over
  $k = \complexnums$. Working over more general base fields does not make a big difference for our arguments
  due to lemma \ref{dense} below, but it has some minor technical advantages; see for example remark
  \ref{why_det}.ii.
\end{rem}
By an \emph{algebraic stack} $\stackM$ over $k$, we always mean an Artin stack $\stackM$ that is locally of
finite type over $k$ (but not necessarily quasi-compact); a standard reference for these notions is
\cite{laumon}. Recall that $\stackM$ can be given by a groupoid $\stackM( S)$ for each $k$-scheme $S$, a
functor $f^*: \stackM(S) \to \stackM(T)$ for each morphism of $k$-schemes $f: T \to S$ and isomorphisms of
functors $(f \circ g)^* \cong g^* \circ f^*$. Such stacks over $k$ form a $2$-category: A $1$-morphism $\Phi:
\stackM \to \stackM'$ can be given by functors $\Phi( S): \stackM( S) \to \stackM'( S)$ and isomorphisms
of functors $f^* \circ \Phi( S) \cong \Phi( T) \circ f^*$; a $2$-morphism $\tau: \Phi_1 \Rightarrow \Phi_2$
can be given by natural transformations $\tau( S): \Phi_1( S) \Rightarrow \Phi_2( S)$. In all this, one can
replace the base field $k$ more generally by a noetherian ring $A$.

Two irreducible algebraic stacks $\stackM$ and $\stackM'$ over $k$ are called \emph{birational} if some open
substacks $\emptyset \neq \stackU \subseteq \stackM$ and $\emptyset \neq \stackU' \subseteq \stackM'$ are
$1$-isomorphic, and \emph{stably birational} if $\Aa^s \times \stackM$ and $\Aa^{s'} \times \stackM'$ are
birational for some $s, s' \in \naturals$.
\begin{example}
  We denote by $\Bunstack_{r, d}$ the moduli stack of vector bundles $E$ on $C$ with rank $r \geq 1$ and
  degree $d \in \integers$. More precisely, $\Bunstack_{r, d}( S)$ is the groupoid of all vector bundles
  $\stackE$ on $C \times_k S$ with rank $r$ and degree $d$ over every geometric point of $S$, and $f^*:
  \Bunstack_{r, d}( S) \to \Bunstack_{r, d}( T)$ is `the' usual pull-back functor.
\end{example}

\begin{example}
  Given a line bundle $L$ of degree $d$ on $C$, we denote by $\Bunstack_{r, L} \subseteq \Bunstack_{r, d}$
  the (reduced) closed substack of vector bundles $E$ on $C$ with $\det( E) \cong L$. More precisely,
  $\Bunstack_{r,L}( S) \subseteq \Bunstack_{r,d}( S)$ is the full subgroupoid of all vector bundles $\stackE$
  on $C \times_k S$ for which $\det( \stackE)$ is Zariski-locally in $S$ isomorphic to the pull-back of $L$.
\end{example}

The basic properties of the algebraic stacks $\Bunstack_{r, d}$ and $\Bunstack_{r, L}$ can be found for
example in \cite{faltings}, \cite{jochen} or the appendix of \cite{rat}: They are smooth over $k$ by standard
deformation theory, $\dim_k \Bunstack_{r, d} = r^2( g-1)$ and $\dim_k \Bunstack_{r, L} = r^2( g-1) - g$.
Moreover, they are geometrically connected, hence geometrically irreducible.
\begin{rems} \label{why_det}
  i) Sending a vector bundle $E$ on $C$ to the line bundle $\det( E)$ defines a smooth morphism $\det:
  \Bunstack_{r, d} \to \Pic^d( C)$ to the Picard \emph{scheme} of $C$. Its fibre over the $k$-point given by
  a line bundle $L$ of degree $d$ on $C$ is precisely $\Bunstack_{r, L}$.

  ii) Let $K = k( \Pic^d( C))$ be the function field, and let the line bundle $\stackL$ on $C \otimes_k K$ be
  the generic fibre of a Poincar\'{e} family. Then the moduli stack $\Bunstack_{r, \stackL}$ over $K$ is
  precisely the generic fibre of $\det$. Thus birational results for fixed degree moduli spaces follow from
  their analogues for fixed determinant moduli spaces; this allows us to concentrate on the latter in the
  sequel.
\end{rems}

\begin{lemma} \label{dense}
  Every open substack $\emptyset \neq \stackU \subseteq \Bunstack_{r, L}$ contains a $k$-rational point.
\end{lemma}
\begin{proof}
  This follows from the well-known unirationality of $\Bunstack_{r, L}$. More precisely:

  We consider the open substack in $\Bunstack_{r, L}$ of vector bundles $E$ with $\cohom^1( E) = 0$ which are
  generated by $\cohom^0( E)$. Twisting by a sufficiently ample line bundle on $C$, we may assume that it is
  nonempty; we may even assume that it contains $\stackU$ by shrinking the latter.

  Now let $U \subseteq \Hom( L^{\dual}, \O_C^{r+1})$ be the open subscheme of nonzero morphisms $\varphi:
  L^{\dual} \to \O_C^{r+1}$ whose cokernel is torsionfree. Sending $\varphi$ to its cokernel defines a
  morphism $\Phi: U \to \Bunstack_{r, L}$; we claim that its image contains $\stackU$. This implies that
  $\Phi^{-1}( \stackU)$ is nonempty and open in the affine space $\Hom( L^{\dual}, \O_C^{r+1})$ over $k$, so
  it contains a $k$-rational point; its image under $\Phi$ is then a $k$-rational point in $\stackU$.

  To prove the claim, we may assume that $k$ is algebraically closed. Let $E$ be a rank $r$ vector bundle on
  $C$ which is generated by $\cohom^0( E)$; it suffices to construct a surjection $\O_C^{r+1}
  \twoheadrightarrow E$ since its kernel will be isomorphic to $\det( E)^{\dual}$. We follow an argument of
  Ramanan \cite[Lemma 3.1]{ramanan} who attributes the result to Atiyah.

  Let $\Delta \subseteq C \times \cohom^0( E)^{r+1}$ be the closed subscheme where the $r+1$ universal
  sections do not generate $E$. The restriction $\Delta_x$ of $\Delta$ to a point $x \in C( k)$ is the
  inverse image under the surjective evaluation map $\cohom^0( E)^{r+1} \to E_x^{r+1}$ of all $r+1$-tuples
  that do not generate the fibre $E_x \cong k^r$; by linear algebra, the latter has codimension $2$ in
  $E_x^{r+1}$, so $\Delta_x$ has codimension $2$ in $\cohom^0( E)^{r+1}$. Thus $\Delta$ has codimension $2$
  in $C \times \cohom^0( E)^{r+1}$, so its projection to $\cohom^0( E)^{r+1}$ has codimension $\geq 1$.
  This shows that there is indeed a surjection $\O_C^{r+1} \twoheadrightarrow E$.
\end{proof}

\begin{rem}
  The above lemma implies that $k$-points are dense in $\Bunstack_{r, d}$ if and only if they are dense in
  the Jacobian $J_C = \Pic^0( C)$. The latter is easily seen to hold if $C( k)$ is infinite. On the other
  hand, an example with $J_C( k)$ finite due to Faddeev \cite{faddeev} is the Fermat curve
  $C \subseteq \mathbb{P}^2$ with equation $x^5 + y^5 = z^5$ over $k = \rationals$.
\end{rem}

\begin{rem}
  The same arguments as for $\Bunstack_{r, d}$ show in fact that the moduli stack $\Cohstack_{r, d}$ of
  coherent sheaves $F$ on $C$ with rank $r$ and degree $d$ is algebraic, smooth of dimension $(g-1)r^2$ over
  $k$ and irreducible; cf. for example the appendix of \cite{rat}. One still has a determinant morphism
  $\det: \Cohstack_{r, d} \to \Pic^d( C)$, e.\,g.\ by \cite[p. 37]{huybrechts-lehn}.

  This includes in particular the case $r = 0$ and $d \geq 0$, whereas of course $\Cohstack_{0, d} =
  \emptyset$ for $d < 0$.  If $r \geq 1$, then $\emptyset \neq \Bunstack_{r, d} \subseteq \Cohstack_{r,d}$
  is an open substack, thus dense.
\end{rem}

The following further examples of algebraic stacks will also be used in the sequel.
\begin{example}
  We denote by $\Vect_n$ the moduli stack of $n$-dimensional vector spaces over $k$; more precisely,
  $\Vect_n( S)$ is the groupoid of rank $n$ vector bundles on $S$. $\Vect_n$ is $1$-isomorphic to the stack
  quotient $\BGL_n$ of $\Spec( k)$ modulo the trivial action of $\GL_n$; thus it is algebraic, smooth of
  dimension $-n^2$ over $k$ and irreducible.
\end{example}

\begin{example}
  Let $A$ be a finitely generated $k$-algebra, and let $F_{\nu}$ be for $\nu = 1, 2$ a coherent sheaf on $C
  \times_k \Spec(A)$ which is flat over $A$ and whose restriction to each point of $\Spec(A)$ has rank $r_{
  \nu}$ and degree $d_{\nu}$. We denote by $\Extstack(F_2 ,F_1)$ the moduli stack over $A$ of extensions of
  $F_2$ by $F_1$; more precisely, $\Extstack(F_2, F_1)(S)$ is the following groupoid for each $A$-scheme
  $p:S \to \Spec(A)$:
  \begin{itemize}
   \item Its objects are the exact sequence $0 \to p^* F_1 \to \stackF \to p^* F_2 \to 0$ of coherent sheaves
    on $C \times_k S$; note that $\stackF$ is automatically flat over $S$.
   \item Its morphisms are the $\O_{C \times_k S}$-module isomorphisms $\stackF \to \stackF'$ that
    are compatible with the identity on $p^* F_{\nu}$ for $\nu = 1, 2$.
  \end{itemize}
\end{example}
\begin{lemma} \label{Extstack}
  $\Extstack( F_2, F_1)$ is algebraic, of finite type and smooth of relative dimension
  $(g-1) r_1 r_2 + r_1 d_2 - r_2 d_1$ over $A$, with all its fibers geometrically irreducible.
\end{lemma}
\begin{proof}
  Using \cite[Prop. 2.1.10]{huybrechts-lehn}, we can find an injective morphism $i: E' \hookrightarrow E$ of
  vector bundles on $C\times_k \Spec(A)$ with cokernel $F_2$. One checks immediately that
  \begin{equation} \label{Ext_cart} \xymatrix{
    \Extstack( F_2, F_1) \ar[r] \ar[d] & \Extstack( E,  F_1) \ar[d]^{i^*}\\
    \Spec( A) \ar[r]^-0                & \Extstack( E', F_1)
  } \end{equation}
  is a $1$-cartesian diagram; here $0$ is given by the trivial extension $F_1 \oplus E'$.

  We choose an open affine covering $X = U \cup V$ and put $U_A := U \times_k \Spec( A)$, $V_A := V \times_k
  \Spec( A)$. This yields a \v{C}ech complex
  \begin{equation*}
    C^{\bullet}( E, F_1) = \Hom( E|_{U_A}, F_1|_{U_A}) \oplus \Hom( E|_{V_A}, F_1|_{V_A}) \longto[\delta]
                           \Hom( E|_{U_A \cap V_A}, F_1|_{U_A \cap V_A})
  \end{equation*}
  whose differential $\delta$ sends two morphisms $\alpha_U$ over $U_A$ and $\alpha_V$ over $V_A$ to
  $\alpha_U-\alpha_V$ over $U_A \cap V_A$. We denote the associated Picard stack \cite[Exp. XVIII, 1.4]{sga4}
  by $C^1( E, F_1)/C^0( E, F_1)$; this is by definition the stack generated by the prestack $\stackM$
  consisting of the following groupoid $\stackM( S)$ for each $A$-scheme $p: S \to \Spec( A)$:
  \begin{itemize}
   \item Its objects are the $1$-cochains $\varphi \in C^1( E, F_1) \otimes_A \Gamma( S, \O_S)$.
   \item Its morphisms $\alpha: \varphi \to \psi$ are the $0$-cochains $\alpha \in C^0( E, F_1) \otimes_A
    \Gamma( S, \O_S)$ with $\delta( \alpha) = \psi - \varphi$.
   \item Its composition law is the addition of $0$-cochains.
  \end{itemize}

  If an object $\varphi$ of $\stackM( S)$ is given, then $\id + \varphi$ is an automorphism of the trivial
  extension $p^*F_1 \oplus p^*E$ over $(U \cap V) \times_k S$; using it to glue the trivial extensions over
  $U \times_k S$ and over $V \times_k S$, we obtain an object of $\Extstack(E, F_1)(S)$ which we denote by
  $\stackF_{\varphi}$. Sending $\varphi$ to $\stackF_{\varphi}$ defines a fully faithful functor
  \begin{equation} \label{glue_S}
    \stackM(S) \longto \Extstack(E, F_1)(S)
  \end{equation}
  because for every morphism $\alpha: \stackF_{\varphi} \to \stackF_{\psi}$ in $\Extstack( E, F_1)(S)$, its
  restrictions $\alpha_U$ to $U \times_k S$ and $\alpha_V$ to $V \times_k S$ form a morphism $\varphi \to
  \psi$ in $\stackM( S)$ and conversely.

  The functors \eqref{glue_S} are compatible with pull-backs along morphisms $f: T \to S$ of schemes over
  $A$, so they induce a $1$-morphism of stacks over $A$
  \begin{equation} \label{glue}
    C^1( E, F_1)/C^0( E, F_1) \longto \Extstack( E, F_1)
  \end{equation}
  which is in fact an isomorphism because \eqref{glue_S} is essentially surjective for affine $S$, since
  $p^* E$ is then projective over $U \times_k S$ and over $V \times_k S$. The same holds for $E'$ instead of
  $E$, and the $1$-isomorphism \eqref{glue} commutes with $i^*$ by construction.

  We denote the mapping cone of $i^*: C^{\bullet}( E, F_1) \to C^{\bullet}( E', F_1)$ by $C^{\bullet}( F_2,
  F_1)$; this is a length $2$ complex of flat $A$-modules whose cohomology is finitely generated over $A$ by
  EGA III \cite[2.2.1]{ega3}. Hence there is a quasi-isomorphism
  \begin{equation*}
    V^{\bullet}=\big[V^0 \longto[\delta^0] V^1 \longto[\delta^1] V^2 \big] \longto[\sim] C^{\bullet}(F_2,F_1)
  \end{equation*}
  from a length $2$ complex of finitely generated flat $A$-modules $V^{\nu}$, i.\,e.\ of vector bundles
  on $\Spec( A)$; cf. \cite[6.10.5 and 7.7.12.i]{ega3}. For each $A$-algebra $B$, we have $\cohom^{\nu}(V^{
  \bullet} \otimes_A B) \cong \Ext^{\nu}(F_2 \otimes_A B, F_1 \otimes_A B)$ by construction. If $B$ is a
  field, then $\Ext^2$ vanishes since $C$ is a smooth curve; hence $\delta^1: V^1 \to V^2$ is surjective by
  Nakayama's lemma, and replacing $V^1$ by its kernel if necessary, we may assume $V^2 = 0$. In this
  situation, it is straightforward to check that the diagram of Picard stacks
  \begin{equation} \label{chains_cart} \xymatrix{
    V^1/V^0 \ar[r] \ar[d] & C^1( E,  F_1)/C^0( E,  F_1) \ar[d]^{i^*}\\
    \Spec( A) \ar[r]^-{0} & C^1( E', F_1)/C^0( E', F_1)
  } \end{equation}
  is $1$-cartesian as well; here $V^1/V^0$ is by definition the stack quotient of the total space $\Tot(V^1)$
  modulo the action of the algebraic group $V^0/A$ given by $\delta^0$. Comparing the diagrams
  \eqref{Ext_cart} and \eqref{chains_cart} shows that $\Extstack(F_2, F_1)$ and $V^1/V^0$ are $1$-isomorphic;
  using the Riemann-Roch formula for $\rank( V^0) - \rank( V^1)$, the lemma follows.
\end{proof}

By a \emph{vector bundle} $\stackV$ of rank $n$ on a stack $\stackM$ over $k$, we mean a $1$-morphism
$\stackV: \stackM \to \Vect_n$. So $\stackV$ assigns to each object $\stackE$ of $\stackM( S)$ a vector
bundle $\stackV( \stackE)$ on the $k$-scheme $S$, in a suitably functorial way; cf.
\cite[Prop. 13.3.6]{laumon}. More explicitly, $\stackV$ is given by a functor $\stackV( S): \stackM( S) \to
\Vect_n( S)$ for each $k$-scheme $S$, and an isomorphism of functors $f^* \circ \stackV( S) \cong \stackV( T)
\circ f^*$ for each $k$-morphism $f: T \to S$.
\begin{examples} \label{Hombdl}
  We fix $r_1, r_2 \geq 1$ and two line bundles $L_1, L_2$ on $C$.

  i) By semicontinuity, there is an open substack $\emptyset \neq \stackU \subseteq \Bunstack_{r_1, L_1}
  \times \Bunstack_{r_2, L_2}$ of vector bundles $E_1, E_2$ for which $\dim \Ext^1(E_1,E_2)$ is minimal, say
  equal to $e$. The vector spaces $\Hom(E_1,E_2)$ are the fibers of a vector bundle $\Hom(\Euniv_1,\Euniv_2)$
  on $\stackU$ due to the theory of cohomology and base change in EGA III \cite[Prop. 7.8.4]{ega3}, using
  that $\stackU$ is locally noetherian and reduced. According to Riemann-Roch, its rank is
  $\rank( \Hom( \Euniv_1, \Euniv_2)) = (1-g) r_1 r_2 + r_1 \deg( L_2) - r_2 \deg( L_1) + e$.

  ii) We moreover fix a vector bundle $E_1$ on $C$ with $\rank( E_1) = r_1$ and $\det( E_1) \cong L_1$.
  Similarly, there is an open substack $\emptyset \neq \stackU_{E_1}' \subseteq \Bunstack_{r_2, L_2}$
  where $\dim \Ext^1( E_1, E_2)$ is minimal, say equal to $e_{E_1}$, and on that a vector bundle $\Hom( E_1,
  \Euniv_2)$ with fibers $\Hom(E_1, E_2)$. This $\Hom(E_1, \Euniv_2)$ is the restriction of $\Hom( \Euniv_1,
  \Euniv_2)$ to the moduli point of $E_1$ if the latter is general, more precisely if $e_{E_1} = e$.

  iii) Arguing as in ii for the other variable, we also obtain for each vector bundle $E_2$ on $C$ a vector
  bundle $\Hom( \Euniv_1, E_2)$ on an open substack $\emptyset \neq \stackU_{E_2}'' \subseteq \Bunstack_{r_1,
  L_1}$.
\end{examples}

As in \cite{king-schofield}, a key ingredient to our main result will be a variant of Hirschowitz' theorem
\cite{hirschowitz} that the tensor product of two general vector bundles is nonspecial. More precisely, we
will use the following:
\begin{prop} \label{hirschowitz}
  Given $r_1 > r_2 \geq 1$ and two line bundles $L_1$, $L_2$ on $C$ with
  \begin{equation} \label{slopes}
    \deg(L_2)/r_2 - \deg(L_1)/r_1 > g-1,
  \end{equation}
  there is a surjective morphism $\varphi: E_1 \twoheadrightarrow E_2$ of vector bundles on $C$ with $\rank(
  E_{\nu}) = r_{\nu}$ and $\det( E_{\nu}) \cong L_{\nu}$ such that $\Ext^1( E_1, E_2) = 0$.
\end{prop}
\begin{proof}
  Like \cite[Lemma 2.1]{king-schofield}, this proof is similar to the one given by Russo and Teixidor
  \cite[Thm. 1.2]{russo-teixidor} for Hirschowitz' theorem itself. Since these references don't cover the
  fixed determinant case that we need, we recall the main arguments here.

  Let $\emptyset \neq \stackU \subseteq \Bunstack_{r_1, L_1} \times \Bunstack_{r_2,L_2}$ be the open substack
  of vector bundles $E_1$ and $E_2$ with $\dim \Ext^1(E_1,E_2)$ minimal, say equal to $e$. According to
  generic flatness \cite[\S 6.9]{ega4}, there is an open substack $\emptyset \neq \stackM$ in the total space
  of the vector bundle $\Hom( \Euniv_1, \Euniv_2)$ over $\stackU$ such that the cokernel of the universal
  family $\varphi^{\univ}$ of morphisms $\varphi: E_1 \to E_2$ is flat over $\stackM$. Then the image of
  $\varphi^{\univ}$ is is also flat over $\stackM$, hence a vector bundle over $C \times_k \stackM$, say of
  rank $r$ and of degree $d$ over every point of $\stackM$. Since $\stackM( k) \neq \emptyset$ according to
  lemma \ref{dense}, it suffices to show $e = 0$, $r = r_2$ and $d = \deg( L_2)$. Note that $r > 0$ by
  assumption \eqref{slopes} and Riemann-Roch.

  By construction, $\stackM$ is an irreducible smooth algebraic stack over $k$ with
  \begin{equation} \label{absdim}
    \dim \stackM = (g-1)( r_1^2 + r_2^2) - 2g + (1-g) r_1 r_2 + r_1 \deg( L_2) - r_2 \deg( L_1) + e
  \end{equation}
  due to Riemann-Roch. On the other hand, sending $\varphi: E_1 \to E_2$ to its kernel, image and cokernel
  defines a $1$-morphism
  \begin{equation*}
    \stackM \longto \stackD \subseteq \Bunstack_{r_1 - r, \deg( L_1) - d} \times \Bunstack_{r, d}
                              \times \Cohstack_{r_2 - r, \deg( L_2) - d}
  \end{equation*}
  to the closed substack $\stackD$ of triples $F_1$, $F$, $F_2$ with $\det( F_{\nu}) \otimes \det( F) \cong
  L_{\nu}$ for $\nu = 1, 2$. Now $\stackD$ is irreducible of codimension $2g$ (since $r_1 - r > 0$ and
  $r > 0$), so
  \begin{equation} \label{basedim}
    \dim \stackD = (g-1) \big[ (r_1 - r)^2 + r^2 + (r_2 - r)^2 \big] - 2g.
  \end{equation}
  The points $\varphi: E_1 \to E_2$ in $\stackM$ with given image in $\stackD$ correspond to extensions
  \begin{equation*}
    0 \to \ker( \varphi) \to E_1 \to    \im( \varphi) \to 0 \quad\text{and}\quad
    0 \to  \im( \varphi) \to E_2 \to \coker( \varphi) \to 0
  \end{equation*}
  that satisfy the open conditions `$E_2$ torsionfree' and `$\dim \Ext^1( E_1, E_2) \leq e$'; using lemma
  \ref{Extstack}, this means that $\stackM$ is smooth over $\stackD$ of relative dimension
  \begin{equation} \label{reldim} \begin{split}
    \dim \stackM - \dim \stackD & = \quad (g-1)(r_1 - r) r + (r_1 - r) d - r (\deg L_1 - d)\\
                                & \quad + (g-1)r (r_2 - r) + r (\deg L_2 - d) - (r_2 - r) d.
  \end{split} \end{equation}
  Combining the dimension formulas \eqref{absdim}, \eqref{basedim} and \eqref{reldim} yields
  \begin{equation} \label{cmpdim}
    (r_1 - r)( \deg L_2 - d) - (r_2 - r)( \deg L_1 - d) + e = (g-1)(r_1 - r)(r_2 - r).
  \end{equation}

  Now suppose $r < r_2$. The image of $\stackM$ in $\Bunstack_{r_1, L_1} \times \Bunstack_{r_2, L_2}$ is open
  and dense, so it contains a pair $E_1$, $E_2$ of \emph{stable} vector bundles on $C$. Consequently, the
  slopes $\mu_{\nu} := (\deg L_{\nu} - d)/(r_{\nu} - r)$ of $\ker( \varphi)$ and $\coker( \varphi)$ satisfy
  $\mu_1 < \deg( L_1)/r_1$ and $\deg( L_2)/r_2 < \mu_2$. On the other hand, $\mu_2 - \mu_1 \leq g-1$ by
  equation \eqref{cmpdim}. This leads to a contradiction with our assumption \eqref{slopes}, thereby proving
  $r=r_2$.

  Now equation \eqref{cmpdim} reads $(r_1 - r)( \deg L_2 - d) + e = 0$. But $d \leq \deg( L_2)$ because
  otherwise $\Cohstack_{r_2 - r, \deg( L_2) - d}$ would be empty, $r_1 - r \geq r_1 - r_2 > 0$ by hypothesis
  and of course also $e \geq 0$; thus we conclude $\deg( L_2) = d$ and $e = 0$.
\end{proof}

\section{$\Gm$-stacks and $\Gm$-gerbes} \label{gerbes}

\begin{defn}
  A \emph{$\Gm$-stack} $\stackM = (\stackM,\iota)$ over $k$ consists of an algebraic stack $\stackM$ over $k$
  together with a group homomorphism $\iota( \stackE): \Gamma( S, \O_S^*) \to \Aut_{\stackM( S)}( \stackE)$
  for each $k$-scheme $S$ and each object $\stackE$ of the groupoid $\stackM( S)$ such that the diagrams
  \begin{equation*} \xymatrix{
      \Gamma( S, \O_S^*) \ar[r]^-{\iota( \stackE)} \ar[dr]_{\iota( \stackE')} & \Aut_{\stackM( S)}( \stackE)
      \ar[d]^{\alpha \mapsto \varphi \alpha \varphi^{-1}}\\ & \Aut_{\stackM( S)}( \stackE')
    } \quad\text{and}\quad \xymatrix{
      \Gamma( S, \O_S^*) \ar[r]^-{\iota( \stackE)} \ar[d]_{f^*} & \Aut_{\stackM( S)}( \stackE) \ar[d]^{f^*}\\
      \Gamma( T, \O_T^*) \ar[r]^-{\iota( f^* \stackE)}          & \Aut_{\stackM( T)}( f^* \stackE)
  } \end{equation*}
  commute for each morphism $\varphi:\stackE \to \stackE'$ in $\stackM(S)$ and each $k$-morphism $f:T \to S$.
\end{defn}

\begin{examples}
  The moduli stacks of vector bundles $\Bunstack_{r,d}$ and $\Bunstack_{r,L}$ are in a canonical way
  $\Gm$-stacks: The multiplication by scalars defines a canonical homomorphism $\iota( \stackE): \Gamma( S,
  \O_S^*) \to \Aut( \stackE)$ for each $k$-scheme $S$ and each vector bundle $\stackE$ on $C \times_k S$;
  the functoriality conditions above are obviously satisfied.
\end{examples}

\begin{defn}
  A $\Gm$-stack $(\stackM, \iota)$ is a \emph{$\Gm$-gerbe} if $\iota( \stackE): \Gamma( S, \O_S^*) \to \Aut_{
  \stackM(S)}(\stackE)$ is an isomorphism for every $k$-scheme $S$ and every object $\stackE$ of $\stackM(
  S)$.
\end{defn}

\begin{examples}
  The open substacks $\Bunstack^{\stab}_{r, d} \subseteq \Bunstack_{r, d}$ and $\Bunstack^{\stab}_{r, L}
  \subseteq \Bunstack_{r, L}$ of geometrically stable vector bundles $E$ on $C$ are $\Gm$-gerbes.
\end{examples}
\begin{proof}
  Geometrically stable vector bundles $E$ in $C$ are known to be \emph{simple}, i.\,e.\ $\End( E) = k \cdot
  \id_E$. If $S$ is a $k$-scheme and $\stackE$ is a vector bundle on $C \times_k S$ whose
  restriction to each point of $S$ is simple, then the theory of cohomology and base change in EGA III
  \cite{ega3} shows $\End( \stackE) = \Gamma( S, \O_S) \cdot \id_{ \stackE}$. (Strictly speaking, \cite{ega3}
  applies only if $S$ is locally noetherian, but standard limit arguments allow us to assume that.) Hence
  $\iota( \stackE): \Gamma( S, \O_S^*) \to \Aut( \stackE)$ is indeed an isomorphism.
\end{proof}

\begin{example}
  If $\coarseM$ is an algebraic space locally of finite type over $k$, then the stack $\stackM := \coarseM
  \times_k \BGm$ is a $\Gm$-gerbe over $k$ in an obvious way.
\end{example}

\begin{defn}
  Let $(\stackM, \iota)$ and $(\stackM', \iota')$ be $\Gm$-stacks over $k$. A $1$-morphism $\Phi: \stackM \to
  \stackM'$ \emph{has weight $w \in \integers$} if the diagram
  \begin{equation*} \xymatrix{
    \Gamma(S,\O_S^*) \ar[r]^-{\iota(\stackE)}\ar[d]_{(\_)^w} &\Aut_{\stackM(S)}(\stackE)\ar[d]^{\Phi(S)}\\
    \Gamma(S,\O_S^*) \ar[r]^-{\iota'(\Phi(\stackE))} &\Aut_{\stackM'( S)}( \Phi( \stackE))
  } \end{equation*}
  commutes for every $k$-scheme $S$ and every object $\stackE$ of the groupoid $\stackM( S)$.
\end{defn}

\begin{example} \label{weight-1}
  There is a canonical $1$-isomorphism $\Bunstack_{r, L} \to \Bunstack_{r, L^{\dual}}$ of stacks over $k$
  which has weight $-1$: It sends a vector bundle $\stackE$ on $C \times_k S$ to $\stackE^{\dual}$ and
  an isomorphism $\varphi: \stackE_1 \to \stackE_2$ of vector bundles on $C \times_k S$ to
  $(\varphi^{-1})^{\dual}: \stackE_1^{\dual} \to \stackE_2^{\dual}$.
\end{example}

\begin{lemma} \label{weight+-1}
  Let $(\stackM, \iota)$ and $(\stackM', \iota')$ be $\Gm$-gerbes over $k$. If $\stackM$ is connected, then
  every $1$-isomorphism $\Phi: \stackM \to \stackM'$ has weight $1$ or has weight $-1$.
\end{lemma}
\begin{proof}
  Let $S$ be a $k$-scheme, and let $\stackE$ be an object of $\stackM( S)$. For each morphism of $k$-schemes
  $f: T \to S$, we consider the group isomorphism
  \begin{equation*}
    \Gamma( T, \O_T^*) \cong \Aut_{\stackM(  T)}(       f^* \stackE) \xrightarrow{\Phi( T)}
                             \Aut_{\stackM'( T)}( \Phi( f^* \stackE)) \cong \Gamma( T, \O_T).
  \end{equation*}
  As $T$ varies, these form an automorphism of the group scheme $\Gm$ over $S$. Hence $S$ is a disjoint union
  of open subschemes $S_+$ and $S_-$ where this automorphism is the identity and $\lambda \mapsto \lambda^{-1
  }$, respectively. This shows that $\stackM$ is the disjoint union of open substacks $\stackM_+$ and
  $\stackM_-$ where $\Phi$ has weight $1$ and weight $-1$, respectively; since $\stackM$ is connected, one of
  these is empty.
\end{proof}

Recall that the \emph{coarse space} $\coarseM$ associated to an algebraic stack $\stackM$ over $k$ is the
fppf-sheafification of the functor $(\schemes/k)^{\op} \to \sets$ that sends each $k$-scheme $S$ to the set
of isomorphism classes in the groupoid $\stackM( S)$. We say that $\coarseM$ is a scheme or an algebraic
space if this sheafified functor is representable by such.
\begin{examples}
  The coarse spaces associated to the $\Gm$-gerbes $\Bunstack_{r,d}^{\stab}$ and $\Bunstack_{r,L}^{\stab}$
  are the usual quasiprojective coarse moduli schemes $\Buncoarse_{r,d}$ and $\Buncoarse_{r, L}$ of stable
  vector bundles $E$ on $C$. We denote by $\Buncoarse_{r,d} \subseteq \overline{\Buncoarse}_{r,d}$ and by
  $\Buncoarse_{r, L} \subseteq \overline{\Buncoarse}_{r, L}$ the natural compactifications given by
  semistable vector bundles.
\end{examples}

\begin{rem}
  Each $\Gm$-gerbe $\stackM$ in the sense above is a gerbe with band $\Gm$ over its coarse space $\coarseM$
  in the sense of \cite[D\'{e}f. IV.2.2.2]{giraud}; cf. also \cite[Rem. 3.19]{laumon}. This justifies
  the terminology `$\Gm$-gerbe'.
\end{rem}

\begin{lemma} \label{gerbe}
  Let $\coarseM$ be the coarse space associated to a $\Gm$-gerbe $\stackM$ over $k$.

  i) $\coarseM$ is an algebraic space locally of finite type over $k$.

  ii) The canonical $1$-morphism $\pi: \stackM \to \coarseM$ is faithfully flat of finite presentation.

  iii) The induced map of points $|\pi|: |\stackM| \to |\coarseM|$ is a homeomorphism for the Zariski
  topologies defined by \cite[Def. II.6.9]{knutson} and \cite[5.5]{laumon}.
\end{lemma}
\begin{proof}
  Due to \cite[Cor. 10.8]{laumon}, $\coarseM$ is an algebraic space, and $\pi$ is faithfully flat and
  locally of finite presentation. The induced map $|\pi|: |\stackM| \to |\coarseM|$ is automatically
  continuous, and it is bijective by the definition \cite[D\'{e}f. 5.2]{laumon} of point; $|\pi|$ is also
  open since $\pi$ is flat \cite[Prop. 5.6]{laumon}. This proves iii; in particular, $\pi$ is quasicompact,
  which completes the proof of ii. Using \cite[2.7.1 and 17.7.5]{ega4}, ii implies that $\coarseM$ is locally
  of finite type over $k$ because $\stackM$ is so, by our convention on stacks.
\end{proof}

\begin{rems}
  i) The coarse space $\pi: \stackM \to \coarseM$ associated to a $\Gm$-gerbe $\stackM$ over $k$ is by
  definition universal for morphisms from $\stackM$ to algebraic spaces over $k$; in particular, $\coarseM$
  is reduced if $\stackM$ is. Part iii of the above lemma \ref{gerbe} implies that $\coarseM$ is irreducible
  if $\stackM$ is. Hence $\coarseM$ is integral if $\stackM$ is.

  ii) If irreducible $\Gm$-gerbes $\stackM$ and $\stackM'$ over $k$ are birational, then their coarse spaces
  $\coarseM$ and $\coarseM'$ are birational as well, due to the above lemma \ref{gerbe}.iii again.
\end{rems}

A $\Gm$-gerbe $\stackM$ is called \emph{neutral} if the canonical $1$-morphism to its coarse space $\pi:
\stackM \to \coarseM$ admits a section $s: \coarseM \to \stackM$. If this is the case, then there is a
$1$-isomorphism $\stackM \cong \coarseM \times \BGm$ of weight $1$ according to \cite[Lemme 3.21]{laumon}.

\begin{rem}
  Lemma \ref{gerbe}.ii above implies that every $\Gm$-gerbe $\stackM$ over $k$ is locally neutral for the
  fppf-topology on its coarse space $\coarseM$.
\end{rem}

\begin{defn}
  Let $\pi: \stackM \to \coarseM$ be the coarse space associated to a $\Gm$-gerbe $\stackM$ over $k$, and let
  $p: S \to \coarseM$ be a scheme over $\coarseM$. A \emph{Poincar\'{e} family} on $S$ is an object $\stackE$
  of $\stackM(S)$ whose classifying morphism $c_{\stackE}:S \to \stackM$ satisfies $p=\pi \circ c_{\stackE}$.
\end{defn}
Let $S \to \coarseM \gets \stackM$ still be a scheme over the coarse space $\coarseM$ of a $\Gm$-gerbe
$\stackM$. Then $\stackS := \stackM \times_{\coarseM} S$ is a $\Gm$-gerbe with coarse space $S$; this
$\Gm$-gerbe is neutral if and only if there is a Poincar\'{e} family $\stackE$ on $S$.

\section{Vector bundles on $\Gm$-gerbes} \label{vectorbundles}

\begin{defn}
  Let $( \stackM, \iota)$ be a $\Gm$-stack over $k$. A vector bundle $\stackV$ on $\stackM$ \emph{has weight
  $w \in \integers$} if the diagram
  \begin{equation*} \xymatrix{
    \Gamma(S,\O_S^*) \ar[r]^-{\iota(\stackE)}\ar[d]_{(\_)^w} &\Aut_{\stackM(S)}(\stackE)\ar[d]^{\stackV(S)}\\
    \Gamma(S,\O_S^*) \ar[r]^-{\cdot \id_{\stackV( \stackE)}} &\Aut_{\O_S}( \stackV( \stackE))
  } \end{equation*}
  commutes for every $k$-scheme $S$ and every object $\stackE$ of the groupoid $\stackM( S)$.
\end{defn}

\begin{example}
  The trivial vector bundle $\O^n$ on any $\Gm$-stack $\stackM$ has weight $0$.
\end{example}

\begin{examples}
  Given a point $P \in C( k)$, we denote by $\Euniv_P$ the restriction of the universal vector bundle
  $\Euniv$ on $C \times \Bunstack_{r, L}$ to $\{P\} \times \Bunstack_{r, L} \cong \Bunstack_{r, L}$.
  \begin{itemize}
   \item[i)] $\Euniv_P$ is a vector bundle of weight $1$ on $\Bunstack_{r, L}$.
   \item[ii)] $(\Euniv_P)^{\dual}$ is a vector bundle of weight $-1$ on $\Bunstack_{r, L}$.
  \end{itemize}
\end{examples}

\begin{examples}
  We fix $r \in \naturals$, a line bundle $L$ and a vector bundle $F$ on $C$. Using the notation of the
  examples \ref{Hombdl}.ii and iii, we have:
  \begin{itemize}
   \item[i)] The vector bundle $\Hom(F, \Euniv)$ on $\stackU_F' \subseteq \Bunstack_{r, L}$ has weight $1$.
   \item[ii)] The vector bundle $\Hom(\Euniv, F)$ on $\stackU_F'' \subseteq \Bunstack_{r,L}$ has weight $-1$.
  \end{itemize}
\end{examples}

\begin{example} \label{det}
  The formalism of determinant line bundles \cite{knudsen-mumford} yields a line bundle $\stackL_{\det}$ on
  $\Bunstack_{r, L}$ whose fibre over any point $[E]$ is $\det \cohom^0( E) \otimes \det^{-1} \cohom^1( E)$.
  This $\stackL_{\det}$ is a line bundle of weight $r( 1-g) + \deg(L)$ on $\Bunstack_{r, L}$ by Riemann-Roch.
\end{example}

We recall that vector bundles $\stackV$ on a fixed algebraic stack $\stackM$ over $k$ form a category; a
morphism $\stackV \to \stackW$ of such vector bundles can be described by a morphism $\stackV( \stackE) \to
\stackW( \stackE)$ of vector bundles on the $k$-scheme $S$ for every object $\stackE$ of $\stackM( S)$.
\begin{lemma} \label{descent}
  Let $\pi: \stackM \to \coarseM$ be the coarse space associated to a $\Gm$-gerbe $\stackM$ over $k$.
  Sending $V$ to $\stackV := \pi^*( V)$ defines an equivalence between the category of vector bundles $V$
  on $\coarseM$ and the category of weight $0$ vector bundles $\stackV$ on $\stackM$.
\end{lemma}
\begin{proof}
  We claim that the inverse equivalence is given by the functor $\pi_*$ from quasicoherent $\O_{\stackM}
  $-modules to quasicoherent $\O_{\coarseM}$-modules \cite[Prop. 13.2.6 (iii)]{laumon}. Thus we have to show
  that $\pi_*$ sends weight $0$ vector bundles $\stackV$ on $\stackM$ to vector bundles $V$ on $\coarseM$ and
  that the two adjunction morphisms $\pi^* \pi_* \stackV \to \stackV$ and $V \to \pi_* \pi^* V$ are
  isomorphisms. All this can be verified locally in the fppf-topology on $\coarseM$, so we may assume that
  the $\Gm$-gerbe $\stackM$ is neutral; this special case is easy to check. 
\end{proof}

\begin{example} \label{coarse_det}
  If $\deg(L) = r(g-1)$, then the determinant line bundle $\stackL_{\det}$ on $\Bunstack_{r,L}$ has weight $0
  $; thus its restriction to $\Bunstack_{r,L}^{\stab}$ descends to $\Buncoarse_{r,L}$ by lemma \ref{descent}.
  In fact, $\stackL_{\det}$ descends to a line bundle $\coarseL_{\det}$ on $\overline{\Buncoarse}_{r, L}$ by
  \cite[Prop. 4.2]{Pic_G}.
\end{example}

\begin{cor} \label{gentriv}
  Let $\stackV$ be a vector bundle of weight $0$ on an irreducible $\Gm$-gerbe $\stackM$ over $k$. Then
  $\stackV$ is trivial on some open substack $\emptyset \neq \stackU \subseteq \stackM$.
\end{cor}
\begin{proof}
  We have $\stackV \cong \pi^* (V)$ for a vector bundle $V$ on $\coarseM$. \cite[Prop. II.6.7]{knutson}
  states that some open subspace $\emptyset \neq \coarseU \subseteq \coarseM$ is a scheme; hence $V$ is
  trivial on some open subscheme $\emptyset \neq U \subseteq \coarseU$, and $\stackV$ is trivial on $\stackU
  := \pi^{-1}( U)$.
\end{proof}

\begin{lemma} \label{neutral}
  Let $\stackM$ be a $\Gm$-gerbe over $k$ with coarse space $\pi: \stackM \to \coarseM$. This $\Gm$-gerbe is
  neutral if and only if there is a line bundle $\stackL$ of weight $1$ on $\stackM$.
\end{lemma}
\begin{proof}
  If $\stackM$ is neutral, say $\stackM = \coarseM \times \BGm$, then pulling back the canonical line bundle
  $L^{\univ}$ of weight $1$ on $\BGm$ yields a line bundle of weight $1$ on $\stackM$.

  Conversely, suppose that there is a line bundle $\stackL$ of weight $1$ on $\stackM$, and let $p: \Tot(
  \stackL)^* \to \stackM$ be the complement of the zero section in its total space. We claim that $p$ is a
  section for $\pi$, i.\,e.\ that $\pi \circ p: \Tot( \stackL)^* \to \coarseM$ is a $1$-isomorphism. This can
  be checked locally in the fppf-topology on $\coarseM$, so we may assume $\stackM = \coarseM \times \BGm$.
  Now $\stackL$ and the pull-back of $L^{\univ}$ differ by tensoring with a weight $0$ line bundle, so they
  are isomorphic locally in the fppf-topology on $\coarseM$ by lemma \ref{descent}; hence we may even assume
  that $\stackL$ is the pull-back of $L^{\univ}$. This special case is easy to check.
\end{proof}

\begin{lemma} \label{brauer}
  Let $\stackM$ be an integral $\Gm$-gerbe over $k$ with coarse space $\stackM \to[\pi] \coarseM$.

  i) There is a vector bundle $\stackV_1$ of weight $1$ on some open substack $\emptyset \neq \stackU_1
  \subseteq \stackM$.

  ii) The generic fibre of the algebra bundle $\Endbdl( \stackV_1)$ on $\stackU_1$ descends to a central
  simple algebra $A$ over the function field $k( \coarseM)$.

  iii) The Brauer class $\psi_{\stackM} := [A] \in \Br( k( \coarseM))$ of the central simple algebra $A$ in
  ii does not depend on the choice of $\stackU_1$ and $\stackV_1$ made in i.

  iv) Any two vector bundles $\stackV$, $\stackV'$ of the same rank $n$ and weight $w$ on open substacks
  $\emptyset \neq \stackU, \stackU' \subseteq \stackM$ are isomorphic on some open substack
  $\emptyset \neq \stackU'' \subseteq \stackU \cap \stackU'$.

  v) There is a vector bundle of given rank $n$ and given weight $w$ on some open substack $\emptyset \neq
  \stackU \subseteq \stackM$ if and only if the index of $w \cdot \psi_{\stackM} \in \Br( k( \coarseM))$
  divides $n$.

  vi) There is a Poincar\'{e} family $\stackE$ on some integral, generically finite scheme $S$ of given
  degree $n$ over $\coarseM$ if and only if the index of $\psi_{\stackM}$ itself divides $n$.
\end{lemma}
\begin{proof}
  Replacing the base field $k$ by the function field $k( \coarseM)$ and $\stackM$ by the $\Gm$-gerbe $\stackM
  \times_{\coarseM} \Spec( k( \coarseM))$ over it, we may assume $\coarseM = \Spec( k)$; then the only open
  substack $\emptyset \neq \stackU \subseteq \stackM$ is $\stackU = \stackM$ itself.

  i) Due to Hilbert's Nullstellensatz, there is a finite field extension $k_1 \supseteq k$ with $\stackM(k_1)
  \neq \emptyset$. Then the $\Gm$-gerbe $\stackM_1 := \stackM \times_k \Spec( k_1)$ is neutral, so there is a
  line bundle $\stackL_1$ of weight $1$ on $\stackM_1$ by lemma \ref{neutral}. Since the first projection
  $\pr: \stackM_1 \to \stackM$ is finite, flat and has weight $1$, we get a vector bundle $\pr_* \stackL_1$
  of weight $1$ on $\stackM$.

  ii) Let $\stackV_1$ be a vector bundle of weight $1$ and rank $n_1$ on $\stackU_1 = \stackM$. The vector
  bundle $\Endbdl( \stackV_1)$ of weight $0$ on $\stackM$ and its multiplication descend to an algebra $A$
  of dimension $n_1^2$ over $k$ by lemma \ref{descent}. With $k_1$, $\pr: \stackM_1 \to \stackM$ and
  $\stackL_1$ as in i,
  \begin{equation*}
    A \otimes_k k_1 \cong \End( \pr^* \stackV_1)
                    \cong \End( \pr^* \stackV_1 \otimes \stackL_1^{\dual})
  \end{equation*}
  is a full matrix algebra over $k_1$ because the weight $0$ vector bundle $\pr^* \stackV_1 \otimes
  \stackL_1^{\dual}$ on $\stackM_1$ is trivial by corollary \ref{gentriv}; this shows that $A$ is central
  simple over $k$. 

  iii -- v) If $\stackV$ is a vector bundle of rank $n$ and weight $w$ on $\stackM$, then the weight $0$
  vector bundle $\Hombdl( \stackV_1^{\otimes w}, \stackV)$ on $\stackM$ and the right action of $\Endbdl(
  \stackV_1)^{\otimes w}$ on it descend to a right module $M$ under the $k$-algebra $A^{\otimes w}$ with
  $\dim_k( M) = n_1^w \cdot n$ by lemma \ref{descent}. This defines a Morita equivalence between vector
  bundles $\stackV$ of weight $w$ on $\stackM$ and finitely generated right $A^{\otimes w}$-modules $M$; its
  inverse sends $M$ to $\stackV := \pi^* M \otimes \stackV_1^{\otimes w}$ where the tensor product is taken
  over $\pi^* A^{\otimes w} \cong \End( \stackV_1^{\otimes w})$.

  In particular, the category of finitely generated right $A$-modules is independent of $\stackV_1$ (up to
  equivalence); this implies iii. According to the theory of central simple algebras, a right $A^{\otimes w}
  $-module $M$ is determined up to isomorphism by $\dim_k( M)$, and there is one with $\dim_k( M) = n_1^w
  \cdot n$ if and only if the index of $[A^{\otimes w}]$ divides $n$; this proves iv and v.

  vi) Here $S = \Spec( K)$ for a finite field extension $K \supseteq k$. There is a Poincar\'{e} family
  $\stackE$ on $S = \Spec( K)$ if and only if the $\Gm$-gerbe $\stackM_K := \stackM \times_k \Spec( K)$ is
  neutral, i.\,e. if and only if its Brauer class $\psi_{\stackM_K} \in \Br( K)$ vanishes; the latter holds
  if and only if $K$ is a splitting field for $\psi_{\stackM}$, because $\psi_{\stackM_K} = \psi_{\stackM}
  \otimes_k K$ in general. By the theory of central simple algebras again, $\psi_{\stackM} \in \Br( k)$ has
  a splitting field of given degree $n$ over $k$ if and only if its index divides $n$.
\end{proof}

\begin{rem}
  If $\stackM = \Bunstack_{r, d}$, then $\psi_{\stackM}$ is precisely the Brauer class $\psi_{r, d}$ used by
  King and Schofield \cite[Def. 3.3]{king-schofield}. The above lemma is based on their arguments.
\end{rem}

\begin{rem}
  Much more information about vector bundles on $\Gm$-gerbes, in particular about \emph{their} moduli spaces,
  can be found in M. Lieblich's work \cite{lieblich, lieblich2}.
\end{rem}

\begin{prop} \label{ramanan}
  Suppose that there is a line bundle $\stackL$ of weight $w \in \integers$ on an open substack $\emptyset
  \neq \stackU \subseteq \Bunstack_{r, L}$ with $\deg( L) = 0$. Then $r$ divides $w$.
\end{prop}
\begin{proof}
  Considering $\stackL \otimes_k \bar{k}$ instead of $\stackL$, we may assume that $k$ is algebraically
  closed. Then $L$ is the $r$-th tensor power of a line bundle on $C$, tensoring with which defines an
  isomorphism $\Bunstack_{r, \O} \cong \Bunstack_{r, L}$. Hence we may also assume $L = \O$.

  We denote by $[\O^r]$ the point on $\Bunstack_{r, \O}$ corresponding to the trivial bundle $\O^r$ on $C$.
  Shrinking $\stackU$ if necessary, we may suppose that $\stackU$ is quasicompact; then we can extend
  $\stackL$ to a neighbourhood of the point $[\O^r]$ as follows:

  Choose a quasicompact open substack $\stackU'$ of $\Bunstack_{r, \O}$ that contains $\stackU$ and
  $[\O^r]$. \cite[Cor. 15.5]{laumon} allows us to extend $\stackL$ to a coherent sheaf on $\stackU'$, more
  precisely to a coherent subsheaf $\stackF \subseteq j_* \stackL$ where $j: \stackU \hookrightarrow \stackU'
  $ is the open embedding. Then the double dual $\stackL' := \stackF^{**}$ is a reflexive coherent sheaf of
  rank one on the smooth stack $\stackU'$ and hence a line bundle, cf. \cite[Chap. VII, \S4.2]{bourbaki_AC}.

  Now $\stackL'$ is a line bundle of weight $w$ because its restriction $\stackL$ to $\stackU$ is so. In
  particular, the scalars $\Gm \subseteq \GL_r = \Aut( \O^r)$ act with weight $w$ on the fibre of $\stackL'$
  over $[\O^r]$. But any one-dimensional representation of $\GL_r$ factors through the determinant, so its
  weight $w$ is always a multiple of $r$.
\end{proof}

\begin{cor} \label{index_r}
  We fix a line bundle $L$ of degree $0$ on $C$, $r, n \geq 1$ and $w \in \integers$.

  i) There is a vector bundle $\stackV$ of rank $n$ and weight $w$ on some open substack $\emptyset \neq
  \stackU \subseteq \Bunstack_{r, L}$ if and only if $r$ divides $w \cdot n$.

  ii) There is a Poincar\'{e} family $\stackE$ on some integral, generically finite scheme $S$ of degree $n$
  over $\Buncoarse_{r, L}$ if and only if $r$ divides $n$.
\end{cor}
\begin{proof}
  If $\stackV$ is a vector bundle of rank $n$ and weight $w$ on $\stackU$, then $\stackL := \det( \stackV)$
  is a line bundle of weight $w \cdot n$; the `only if' in i thus follows from the previous proposition.
  According to lemma \ref{brauer}.v, this means that $r$ divides $w$ times the index of the Brauer class
  $w \cdot \psi_{\Bunstack_{r, L}} \in \Br( k( \Buncoarse_{r, L}))$.

  On the other hand, the vector bundle $\Euniv_P$ of rank $r$ and weight $1$ on $\Bunstack_{r, L}$ shows that
  the index of $\psi_{\Bunstack_{r,L}}$ divides $r$. Due to \cite[5.4]{artin}, this implies that the index of
  $w \cdot \psi_{\Bunstack_{r,L}}$ divides $r/w$ if $w$ divides $r$; together with the previous paragraph,
  this shows that the index of $w \cdot \psi_{\Bunstack_{r,L}}$ equals $r/w$ if $w$ divides $r$.

  For general $w$, the Brauer classes $\hcf(r, w) \cdot \psi_{\Bunstack_{r, L}}$ and $w \cdot \psi_{
  \Bunstack_{r, L}}$ have the same splitting fields because $w/\hcf(r, w)$ is invertible modulo $r$; hence
  the index $r/\hcf( r, w)$ of the former is also the index of the latter. Using this, i and ii follow from
  lemma \ref{brauer}.v and vi.
\end{proof}

\section{Grassmannian bundles} \label{Grassmannians}
Given a vector bundle $\stackV$ on a stack $\stackM$ over $k$, we denote by $\Gr_j( \stackV) \to \stackM$ the
associated Grassmannian bundle of $j$-dimensional linear subspaces for $j \leq \rank( \stackV)$. Recall that
this is a smooth representable $1$-morphism of stacks over $k$ whose fibre over the point of $\stackM$
corresponding to an object $E$ of $\stackM( k)$ is the usual Grassmannian scheme $\Gr_j( \stackV( E))$.

If $(\stackM,\iota)$ is a $\Gm$-stack over $k$ and $\stackV$ has some weight $w \in \integers$, then $\Gr_j(
\stackV)$ becomes a $\Gm$-stack over $k$ as well, because all linear subspaces of $\stackV( E)$ are invariant
under the $\Gm$-action given by $\iota$. If $\stackM$ is even a $\Gm$-gerbe over $k$, then $\Gr_j( \stackV)$
also is.
\begin{rem} \label{duality}
  There is a canonical $1$-isomorphism $\Gr_j( \stackV) \longto[\sim] \Gr_{\rank( \stackV)-j}( \stackV^{\dual
  })$ that sends each linear subspace of $\stackV( E)$ to its orthogonal complement in $\stackV( E)^{\dual}$.
\end{rem}

\begin{examples} \label{hecke}
  We denote by $\Parstack_{r, L}^{m/P}$ the moduli stack of rank $r$ vector bundles $E$ on $C$ with $\det( E)
  \cong L$, endowed with a quasiparabolic structure of multiplicity $m$ over the point $P \in C( k)$ in the
  sense of \cite{mehta-seshadri}. We recall that such a quasiparabolic structure can be given by a coherent
  subsheaf $E' \subseteq E$ for which $E/E'$ is isomorphic to the skyscraper sheaf $\O_P^m$.
  \begin{itemize}
   \item[i)] $\Parstack_{r, L}^{m/P}$ is canonically $1$-isomorphic to $\Gr_m((\Euniv_P)^{\dual})$ over
    $\Bunstack_{r, L}$.
   \item[ii)] $\Parstack_{r, L}^{m/P}$ is also $1$-isomorphic to $\Gr_m( \Euniv_P)$ over
    $\Bunstack_{r, L( -mP)}$.
  \end{itemize}
\end{examples}
\begin{proof}
  i) The quasiparabolic vector bundle $E^{\bullet} = (E' \subseteq E)$ is given by the vector bundle
  $E$ together with a dimension $m$ quotient of the fibre $E_P$.

  ii) $E^{\bullet}$ is also given by the vector bundle $E'$ together with a dimension $m$ vector subspace in
  the fibre $E'(P)_P$ of the twisted bundle $E'(P)$. Choosing once and for all a trivialisation at $P$ for
  the line bundle $\O_C( P)$ identifies $E'(P)_P$ and $E_P'$.
\end{proof}

In analogy with \cite{king-schofield}, the proof of our main result will use not only the above Hecke
correspondence, but also the following slightly more involved comparison of Grassmannian bundles over
different moduli stacks $\Bunstack_{r, L}$:
\begin{examples} \label{grass_hecke}
  Given $j, r_1, r_2, r_3 \geq 1$ with $r_1 + r_3 = j r_2$ and line bundles $L_1, L_2, L_3$ on $C$ with $L_1
  \otimes L_3 \cong L_2^{\otimes j}$, let $\stackM$ be the moduli stack of all exact sequences
  \begin{equation} \label{F_seq}
    0 \longto E_1 \longto[i] E_2 \otimes_k V \longto[p] E_3 \longto 0
  \end{equation}
  in which $V$ is a vector space over $k$ and $E_1, E_2, E_3$ are vector bundles on $C$ with $\dim( V) = j$
  and $\rank( E_{\nu}) = r_{\nu}$, $\det( E_{\nu}) \cong L_{\nu}$ for $\nu = 1, 2, 3$. More precisely,
  $\stackM( S)$ is the following groupoid for each $k$-scheme $S$:
  \begin{itemize}
   \item An object consists of a rank $j$ vector bundle $\stackV$ on $S$, an object $\stackE_{\nu}$ of
    $\Bunstack_{r_{\nu}, L_{\nu}}(S)$ for $\nu = 1, 2, 3$ and an exact sequence of bundles on $C \times_k S$
    \begin{equation*}
      0 \longto \stackE_1 \longto \stackE_2 \otimes_{\O_S} \stackV \longto \stackE_3 \longto 0.
    \end{equation*}
   \item A morphism consists of four vector bundle isomorphisms $\lambda: \stackV \to \stackV'$ and
    $\varphi_{\nu}: \stackE_{\nu} \to \stackE_{\nu}'$, $\nu = 1, 2, 3$, such that the following diagram
    commutes:
    \begin{equation*} \xymatrix{
      0 \ar[r] & \stackE_1                          \ar[r] \ar[d]^{\varphi_1}
               & \stackE_2  \otimes_{\O_S} \stackV  \ar[r] \ar[d]^{\varphi_2 \otimes \lambda}
               & \stackE_3                          \ar[r] \ar[d]^{\varphi_3} & 0\\
      0 \ar[r] & \stackE_1'                         \ar[r]
               & \stackE_2' \otimes_{\O_S} \stackV' \ar[r]
               & \stackE_3'                         \ar[r] & 0.
    } \end{equation*}
  \end{itemize}
  $\stackM$ is algebraic over $k$ because the forgetful $1$-morphism $\stackM \to \Bunstack_{r_2, L_2} \times
  \Vect_j$ is representable (by suitable open subschemes of relative $\Quot$-schemes) and of finite type. We
  have a strictly commutative diagram of canonical $1$-morphisms
  \begin{equation*} \xymatrix@C-21ex{
    \Gr_j( \Hom( \Euniv_1, \Euniv_2)) \ar[dr]!<-12ex, 0ex> &&
    \stackM' \subseteq \stackM \supseteq \stackM''
      \ar[ll]_-{\Gamma'} \ar[rr]^-{\Gamma''} \ar[dd]^>>>>>>{\Phi_2}&&
    \Gr_j( \Hom( \Euniv_2, \Euniv_3)) \ar[dl]!<11.5ex, 0ex>\\
    \ar@<1.5ex>@{<-}[urr]+<-7ex, -3ex> & \stackU_{1, 2} \subseteq
    \Bunstack_{r_1, L_1} \times \Bunstack_{r_2, L_2}
      \ar[dr] \ar@{<-}[ur]!< 2ex, 0ex>^<<<<<<{\Phi_{1, 2}} &&
    \Bunstack_{r_2, L_2} \times \Bunstack_{r_3, L_3} \supseteq \stackU_{2, 3}
      \ar[dl] \ar@{<-}[ul]!<-2ex, 0ex>_<<<<<<{\Phi_{2, 3}} &
    \ar@<-1.5ex>@{<-}[ull]+<7ex, -3ex>\\
    && \qquad\qquad\quad \Bunstack_{r_2, L_2} \quad\qquad\qquad &&
  } \end{equation*}
  in which $\Phi_{1, 2}$, $\Phi_2$, $\Phi_{2, 3}$ are forgetful morphisms,
  \begin{itemize}
   \item $\stackU_{\nu,\nu+1} \subseteq \Bunstack_{r_{\nu}, L_{\nu}} \times \Bunstack_{r_{\nu+1}, L_{\nu+1}}$
    is the open substack defined by the condition $\Ext^1( E_{\nu}, E_{\nu+1}) = 0$,
   \item $\stackM'' \subseteq \Phi_{2, 3}^{-1}( \stackU_{2, 3})$ is the open substack where $p_*: V \to \Hom(
    E_2, E_3)$ is injective in \eqref{F_seq}, and $\Gamma''$ sends such a sequence to the image of $p_*$,
   \item $\stackM' \subseteq \Phi_{1,2}^{-1}( \stackU_{1,2})$ is the open substack where $i^*: V^{\dual} \to
    \Hom(E_1, E_2)$ is injective in \eqref{F_seq}, and $\Gamma'$ sends such a sequence to the image of $i^*$,
  \end{itemize}
  and the remaining four unlabelled arrows are canonical projections. We claim:
  \begin{itemize}
   \item[i)] $\Gamma'$ and $\Gamma''$ are open immersions.
   \item[ii)] If $\deg( L_{\nu+1})/r_{\nu+1} - \deg( L_{\nu})/r_{\nu} > g-1$ for $\nu = 1, 2$ and $j$ divides
    $r_{\nu}, \deg( L_{\nu})$ for $\nu = 1, 3$, then $\stackM' \cap \stackM''$ contains a $k$-rational point.
  \end{itemize}
\end{examples}
\begin{proof}
  i) An open substack of $\Gr_j(\Hom( \Euniv_2, \Euniv_3))$ parameterises all linear subspaces $W \subseteq
  \Hom(E_2, E_3)$ for which the evaluation map $\epsilon_W: E_2 \otimes_k W \to E_3$ is surjective.
  $\Gamma''$ is in fact a $1$-isomorphism over $\stackU_{2, 3}$ onto that open substack; its inverse sends
  such a linear subspace $W$ to the sequence
  \begin{equation*}
    0 \longto \ker( \epsilon_W) \longto E_2 \otimes_k W \longto[\epsilon_W] E_3 \longto 0.
  \end{equation*}
  Hence $\Gamma''$ is an open immersion; considering $E_{4-\nu}^{\dual}$ and $V^{\dual}$ instead of $E_{\nu}$
  and $V$, it follows that $\Gamma'$ is an open immersion as well.

  ii) Due to i and lemma \ref{dense}, it suffices to show $\stackM' \cap \stackM'' \neq \emptyset$. Thus we
  may assume that $k = \bar{k}$ is algebraically closed.

  In this case, we can choose a line bundle $\tilde{L}_1$ on $C$ with $\tilde{L}_1^{\otimes j} \cong L_1$;
  we define $\tilde{L}_3 := L_2 \otimes \tilde{L}_1^{\dual}$. Instead of $j$, $r_1, r_2, r_3$ and
  $L_1, L_2, L_3$, we consider $\tilde{j} := 1$, $r_1/j, r_2, r_3/j$ and $\tilde{L}_1, L_2, \tilde{L}_3$.
  We denote the resulting moduli stack by $\tilde{ \stackM}$ and the analogous open substacks by $\tilde{
  \stackM}' \subseteq \tilde{ \stackM} \supseteq \tilde{ \stackM}''$. If $\tilde{ \stackM}' \cap \tilde{
  \stackM}''$ has a $k$-rational point $0 \to \tilde{E}_1 \to E_2 \otimes_k \tilde{V} \to \tilde{E}_3 \to 0$,
  then taking the direct sum of $j$ copies yields a $k$-rational point of $\stackM' \cap \stackM''$.
  Hence it suffices to treat the case $j = 1$.

  In this case $j = 1$, proposition \ref{hirschowitz} implies that the image of $\Gamma''$ described above is
  nonempty; hence $\stackM'' \neq \emptyset$. The dual argument shows $\stackM' \neq \emptyset$. But the
  forgetful $1$-morphism $\stackM \to \Bunstack_{r_1, L_1} \times \Vect_1 \times \Bunstack_{r_3, L_3}$ is
  smooth with all fibers irreducible due to lemma \ref{Extstack}; hence $\stackM$ is irreducible. This
  proves $\stackM' \cap \stackM'' \neq \emptyset$ here.
\end{proof}

\begin{cor} \label{grass_birat}
  In the above notation, let the vector bundle $E_2$ on $C$ belong to a $k$-rational point \eqref{F_seq} of
  $\stackM' \cap \stackM''$. Then the Grassmannian $\Gr_j( \Hom( E_2, \Euniv_3))$ over the open substack
  $\emptyset \neq (\stackU_{2, 3})_{E_2} \subseteq \Bunstack_{r_3, L_3}$ is birational to the Grassmannian
  $\Gr_j( \Hom( \Euniv_1, E_2))$ over the open substack $\emptyset \neq (\stackU_{1, 2})_{E_2} \subseteq
  \Bunstack_{ r_1, L_1}$.
\end{cor}
\begin{proof}
  Both Grassmannian bundles contain as a nonempty open substack the fibre of the forgetful $1$-morphism
  $\Phi_2: \stackM' \cap \stackM'' \to \Bunstack_{r_2, L_2}$ over $E_2$.
\end{proof}

\begin{lemma} \label{rankreduction}
  Suppose that $\stackV$ and $\stackW$ are two vector bundles of the same weight $w$ on an irreducible
  $\Gm$-gerbe $\stackM$ over $k$, and let $j \leq \rank( \stackV)$ be given.

  i) $\Gr_j( \stackV \oplus \stackW)$ is birational to $\Aa^{j \cdot \rank(\stackW)} \times \Gr_j(\stackV)$.

  ii) $\Gr_{j + \rank( \stackW)}( \stackV \oplus \stackW)$ is birational to $\Aa^s \times \Gr_j(\stackV)$
  with $s := (\rank( \stackV) - j) \cdot \rank(\stackW)$.
\end{lemma}
\begin{proof}
  i) Let $\emptyset \neq \stackU \subseteq \Gr_j( \stackV \oplus \stackW)$ be the open substack that
  parameterises all the $j$-dimensional vector subspaces of $\stackV( E) \oplus \stackW( E)$ whose image $S$
  in $\stackV( E)$ still has dimension $j$. Each of these subspaces is the graph of a unique linear map
  $S \to \stackW( E)$. Denoting by $\Suniv \subseteq p^* \stackV$ the universal subbundle on $p: \Gr_j(
  \stackV) \to \stackM$, this defines a canonical $1$-isomorphism between $\stackU$ and the total space of
  the vector bundle $\Hombdl( \Suniv, p^* \stackW)$ on $\Gr_j( \stackV)$. The latter has weight $0$ because
  $\Suniv$ and $p^* \stackW$ both have weight $w$; hence corollary \ref{gentriv} applies to it.

  ii) follows from i by applying remark \ref{duality} twice.
\end{proof}

\begin{prop} \label{weightreduction}
  Let $L$ be a line bundle of degree $0$ on $C$. Let $\stackV$ be a vector bundle of rank $n$ and weight $w$
  on an open substack $\emptyset \neq \stackU \subseteq \Bunstack_{r, L}$, and let $j \leq n$.
  \begin{itemize}
   \item[i)] If $r$ divides $w \cdot j$, then $\Gr_j( \stackV)$ is birational to $\Aa^{j(n-j)} \times
    \Bunstack_{r, L}$.
   \item[ii)] If $w \in \{0, \pm 1\}$, then $\Gr_j( \stackV)$ is birational to $\Aa^s \times \Gr_m(\Euniv_P)$
    for some $s$, where $m \in \{0, 1, \ldots, r-1\}$ is the remainder of $w \cdot j$ modulo $r$.
   \item[iii)] If $r$ and $w \cdot j$ are coprime, then $\Gr_j( \stackV)$ is stably birational to $\Gr_1(
    \Euniv_P)$.
  \end{itemize}
\end{prop}
\begin{proof}
  Note that $r$ divides $w \cdot n$ due to corollary \ref{index_r}.i.

  i) Corollary \ref{index_r}.i also implies that there are weight $w$ vector bundles $\stackV'$ and $\stackW$
  of rank $j$ and $n-j$ on some open substack $\emptyset \neq \stackU' \subseteq \Bunstack_{r, L}$. According
  to lemma \ref{brauer}.iv, we have $\stackV \cong \stackV' \oplus \stackW$ on some nonempty open substack of
  $\stackU \cap \stackU'$; hence $\Gr_j( \stackV)$ is by lemma \ref{rankreduction}.i birational to
  $\Aa^{j(n-j)} \times \Gr_j( \stackV') \cong \Aa^{j(n-j)} \times \stackU'$.

  ii) The case $w=0$ follows from i. Using remark \ref{duality}, it suffices to treat the case $w=1$. Here
  lemma \ref{brauer}.iv allows us to assume $\stackV = (\Euniv_P)^{n/r}$. Then $\Gr_j(\stackV)$ is by lemma
  \ref{rankreduction}.ii birational to $\Aa^{s'} \times \Gr_m( \stackV')$ for some $s'$ with $\stackV' :=
  (\Euniv_P)^{(n-j+m)/r}$. Now the claim follows from lemma \ref{rankreduction}.i.

  iii) Let the algebraic stack $\stackM$ be the fibred product of $p: \Gr_j( \stackV) \to \stackU$ and $q:
  \Gr_1( \Euniv_P) \to \stackU$. Denoting by $\Suniv \subseteq p^* \stackV$ the universal subbundle, $\det(
  \Suniv)$ is a line bundle of weight $w \cdot j$ on $\Gr_j( \stackV)$, whereas $p^* \det( \Euniv_P)$ is a
  line bundle of weight $r$; since these weights are coprime, there also is a line bundle $\stackL$ of weight
  $1$ on $\Gr_j( \stackV)$. Applying corollary \ref{gentriv} to $\stackL^{\dual} \otimes p^* \Euniv_P$, we
  see that the Grassmannian bundle $\stackM \cong \Gr_1( p^* \Euniv_P)$ over $\Gr_j( \stackV)$ is birational
  to $\Aa^{r-1} \times \Gr_j( \stackV)$.

  On the other hand, the universal subbundle of $q^* \Euniv_P$ is a line bundle of weight $1$ on $\Gr_1(
  \Euniv_P)$; applying \ref{gentriv} as above, we see that the Grassmannian bundle $\stackM \cong \Gr_j( q^*
  \stackV)$ over $\Gr_1( \Euniv_P)$ is also birational to $\Aa^{j(n-j)} \times \Gr_1( \Euniv_P)$.
\end{proof}

\begin{rems} \label{weight_x}
  i) Amitsur's conjecture \cite[p. 40]{amitsur} about generic splitting fields of central simple algebras
  would allow to remove `stably' in iii if $j = 1$ and $n = r$. (Indeed, both Grassmannians
  contain \emph{neutral} $\Gm$-gerbes as open substacks here, and their coarse spaces have isomorphic
  function fields according to the conjecture.)

  ii) If $w \notin \{0, \pm 1\}$ and $r$ neither divides nor is coprime to $w \cdot j$, then the index of
  the Brauer class $\psi_{\Gr_j( \stackV)}$ need no longer coincide with its order in the Brauer group,
  according to \cite{schofield-vandenbergh} for $j = 1$ and to \cite{blanchet}, \cite{wadsworth} for general
  $j$. Thus there seems to be no direct generalisation of iii to this case.
\end{rems}

\section{Proof of the main theorem} \label{mainresults}

\begin{thm} \label{mainthm}
  Let $L$ be a line bundle of degree $d$ on $C$. We consider an open substack $\emptyset \neq \stackU_1
  \subseteq \Bunstack_{r,L}$ and a tower of Grassmannian bundles
  \begin{equation*}
    \stackM = \Gr_{j_T   } (\stackV_T    ) \longto \stackU_T \subseteq
              \Gr_{j_{T-1}}(\stackV_{T-1}) \longto \stackU_{T-1} \ldots
              \Gr_{j_1}    (\stackV_1    ) \longto \stackU_1 \subseteq \Bunstack_{r, L}
  \end{equation*}
  where $\stackV_t$ is a vector bundle of weight $w_t \in \{0,\pm 1\}$ and rank $n_t \geq j_t$ on $\stackU_t$
  for each $t$, and $\emptyset \neq \stackU_{t+1} \subseteq \Gr_{j_t}(\stackV_t)$ is again an open substack.

  Let $h$ denote the highest common factor of $r, d, w_1 \cdot j_1, \ldots, w_T \cdot j_T$. Then $\stackM$ is
  birational to $\Aa^s \times \Bunstack_{h, L_0}$ for the degree $0$ line bundle $L_0 := L( -dP)$ on $C$ with
  \begin{equation*} 
    s = \dim \stackM - \dim \Bunstack_{h, L_0} = (g-1)(r^2 - h^2) + \sum_t j_t( n_t -j_t).
  \end{equation*}
\end{thm}
\begin{proof}
  i) We start with the case $h = r$, i.\,e.\ $r$ divides $d$ and $w_1 j_1, \ldots, w_T j_T$. We have a
  $1$-isomorphism $\Bunstack_{r, L_0} \to \Bunstack_{r, L}$ defined by sending each vector bundle $E$ on $C$
  to its twist $E(\frac{d}{r} P)$; hence we may assume $d = 0$ without loss of generality. Then $\Gr_{j_1}(
  \stackV_1)$ is birational to $\Aa^{j_1( n_1 - j_1)} \times \Bunstack_{r, L_0}$ due to proposition
  \ref{weightreduction}.i; replacing $k$ by the function field of $\Aa^{j_1( n_1 - j_1)}$ and iterating,
  the theorem follows here.

  For the remaining case $h < r$, we use induction on $r/h$. We first consider the case that $r$ divides $d$
  and denote by $t \leq T$ the smallest index for which $r$ does not divide $w_t j_t$. Then $\stackU_t$ is
  birational to $\Aa^{s_1} \times \Bunstack_{r, L_0}$ for some $s_1$ by the previous paragraph; replacing $k$
  by the function field of $\Aa^{s_1}$, we may assume $s_1 = 0$ and $t = 1$, i.\,e.\ $\stackU_1 \subseteq
  \Bunstack_{r, L_0}$ and $r$ does not divide $w_1 \cdot j_1$. Here proposition \ref{weightreduction}.ii
  makes $\stackU_2$ birational to $\Aa^{s_2} \times \Gr_m( \Euniv_P)$ for some $s_2$, where $m \in \{1,
  \ldots, r-1\}$ is the remainder of $w_1 \cdot j_1$ modulo $r$. Again, we can achieve $s_2 = 0$ by changing
  $k$. But then $\stackU_2$ is also birational to $\Gr_m( (\Euniv_P)^{\dual}) \to \Bunstack_{r, L_0( mP)}$
  by the examples \ref{hecke}. Thus we may assume for the induction step that $r$ does not divide $d$.

  Let $j$ be the highest common factor of $r$ and $d$; then $r/j > 1$. According to elementary number theory,
  there is a unique integer $r_2$ with $r/j < r_2 < 2r/j$ such that $r_2 d \equiv j \bmod r$; we then define
  $d_2 \in \integers$ by the equation
  \begin{equation} \label{chi=j}
    (1-g) r_2 r + r_2 d - r d_2 = j
  \end{equation}
  and put $L_2:= \O_C(d_2 P)$. Now we claim that the hypotheses of \ref{grass_hecke}.ii are satisfied with
  $r_1:= j r_2 - r < r$, $r_3 := r$ and $L_1 := L^{\dual}( jd_2 P)$, $L_3 := L$:

  Only the two inequalities in \ref{grass_hecke}.ii are not obvious, but $d/r - d_2/r_2 > g-1$ follows from
  \eqref{chi=j}, and then $d_2/r_2 - d_1/r_1 > g-1$ for $d_1 := \deg( L_1)$ follows from the computation
  \begin{equation*}
    r_1 d_2 - r_2 d_1 = (j r_2 - r)d_2 - r_2 (j d_2 - d) = r_2 d - r d_2 > (g-1)r_2 r > (g-1)r_1 r_2.
  \end{equation*}

  Hence corollary \ref{grass_birat} applies for some vector bundle $E_2$ on $C$ of rank $r_2$ and determinant
  $L_2$. It yields that $\Gr_j(\Hom(\Euniv_1,E_2))$ over $(\stackU_{1,2})_{E_2}\subseteq \Bunstack_{r_1,L_1}$
  is birational to $\Gr_j$ of the vector bundle $\Hom( E_2, \Euniv)$ on $(\stackU_{2,3})_{E_2} \subseteq
  \Bunstack_{r,L}$ whose rank is $j$ by \eqref{chi=j}; thus it is birational to $\Bunstack_{r, L}$ itself.
  Example \ref{weight-1} shows that $\Bunstack_{r_1, L_1} \cong \Bunstack_{r_1, L( -jd_2 P)}$; now we can
  apply the induction hypothesis, observing that $r_1 < r$ and $\hcf( r_1, d_1, j) = j = \hcf( r, d)$.
\end{proof}

Keeping the notation of theorem \ref{mainthm} for a while, let $\emptyset \neq \stackU \subseteq \stackM$ be
an open substack which is a $\Gm$-gerbe, and let $\pi: \stackU \to \coarseU$ be its coarse space.
\begin{cor} \label{rationality} \begin{itemize}
  \item[i)] $\coarseU$ is birational to $\Aa^s\times \Buncoarse_{h, L_0}$ over $k$.
  \item[ii)] If $h = 1$, then $\coarseU$ is rational over $k$.
  \item[iii)] If $g = h = 2$ and $\mathrm{char}( k) = 0$, then $\coarseU$ is rational over $k$.
\end{itemize} \end{cor}
\begin{proof}
  i) is a consequence of theorem \ref{mainthm} and lemma \ref{gerbe}.iii.

  ii) follows from i, since $\Buncoarse_{1, L_0} \cong \Spec(k)$.

  iii) Using i, it suffices to show $\overline{\Buncoarse}_{2, L_0} \cong \mathbb{P}^3$, or equivalently
  $\overline{\Buncoarse}_{2, L_0( 2P)} \cong \mathbb{P}^3$. Since $C$, $L$ and $P$ are defined over a
  finitely generated extension field of $\rationals$, we may assume $k \subseteq \complexnums$. Example
  \ref{coarse_det} yields a determinant line bundle $\coarseL_{\det}$ on $\overline{\Buncoarse}_{2,L_0(2P)}$;
  we claim that its dual is very ample and that the resulting projective embedding is an isomorphism onto
  $\mathbb{P}^3$. To check that, we may replace $k$ by its extension field $\complexnums$; then the claim
  holds by \cite[\S 7]{narasimhan-ramanan} and \cite[Prop. 2.5]{beauville}.
\end{proof}

\begin{cor} \label{no_poincare}
  With $\coarseU$ and $h$ as above, the following are equivalent for $n \in \naturals$:
  \begin{itemize}
   \item[i)] There is an integral, generically finite scheme $S$ of degree $n$ over $\coarseU$ that
    admits a Poincar\'{e} family $\stackE$ on $S$.
   \item[ii)] $h$ divides $n$.
  \end{itemize}
  In particular, there is no Poincar{\'e} family on any open subspace of $\coarseU$ if $h \neq 1$.
\end{cor}
\begin{proof}
  Combine theorem \ref{mainthm} with corollary \ref{index_r}.ii.
\end{proof}

\begin{cor} \label{not_birat}
  Let $\stackM$ and $h$ be as in theorem \ref{mainthm}; similarly for $\stackM'$ and $h'$. If $h \neq h'$,
  then the stacks $\stackM$ and $\stackM'$ are not stably birational.
\end{cor}
\begin{proof}
  Suppose that $\Aa^s \times \stackM$ and $\Aa^{s'} \times \stackM'$ have a common open substack $\stackU
  \neq \emptyset$. We may assume that $\stackU$ is a $\Gm$-gerbe; let $\coarseU$ be its coarse space. Due to
  the previous corollary, $h$ and $h'$ both coincide with the minimal degree of integral, generically finite
  schemes $S$ over $\coarseU$ that admit a Poincar\'{e} family; hence $h = h'$.
\end{proof}
\begin{rem}
  Nevertheless, coarse spaces for $\stackM$ and $\stackM'$ might be stably birational for some $h \neq h'$,
  as is the case for $\mathrm{char}(k) = 0$, $g = 2$ and $h = 1$, $h' = 2$ by corollary \ref{rationality}.
  All that \ref{not_birat} says is that an eventual proof of stable birationality for the coarse spaces has
  to contain arguments that do not carry over to the stacks.
\end{rem}

\begin{cor} \label{rank_x_weight}
  Let $\stackM$ and $h$ be as above; let $n \in \naturals$ and $w \in \integers$ be given. Then there exists
  an open substack $\emptyset \neq \stackU \subseteq \stackM$ and a vector bundle $\stackV$ of rank $n$ and
  weight $w$ on $\stackU$ if and only if $h$ divides $w \cdot n$.
\end{cor}
\begin{proof}
  This follows directly from theorem \ref{mainthm}, lemma \ref{weight+-1} and corollary \ref{index_r}.i.
\end{proof}

\begin{example}
  Let $D \subset C( k)$ be a finite set of points on our curve $C$. We consider vector bundles $E$ on $C$
  with fixed rank $r$ and degree $d$, endowed with a quasiparabolic structure with fixed multiplicities
  $m_1(P), \ldots, m_l(P)$ at each point $P \in D$. Let $\Parcoarse_{r, d}^{\underline{m}/D}$ be their coarse
  moduli scheme for some given weights $\alpha_1, \ldots, \alpha_l$, as introduced in \cite{mehta-seshadri}.
  Then the above implies in particular:
  \begin{enumerate}
   \item[i)] All fibers of $\det: \Parcoarse_{r, d}^{\underline{m}/D} \to \Pic^d( C)$ are rational if
    the highest common factor $h$ of $r$, $d$ and all multiplicities $m_l( P)$ is $1$. (This generalises
    \cite{boden-yokogawa}.)
   \item[ii)] There is a Poincar\'{e} family on some nonempty open subscheme of
    $\Parcoarse_{r, d}^{\underline{m}/D}$ if and only if $h = 1$.
  \end{enumerate}
\end{example}

\begin{example} \label{parbdl}
  Concerning the coarse moduli scheme $\Buncoarse_{r, L}$ of vector bundles with fixed determinant $L$ of
  degree $d$, the above reproves the rationality theorem of King and Schofield \cite{king-schofield} and the
  non-existence theorem for Poincar\'{e} families due to Ramanan \cite{ramanan}. Moreover, we get the
  following additional information:
  \begin{itemize}
   \item[i)] There is an integral, generically finite scheme of degree $n$ over $\Buncoarse_{r, L}$ with a
    Poincar\'{e} family if and only if $\hcf(r, d)$ divides $n$.
   \item[ii)] For a fixed vector bundle $F$ on $C$, the minimal dimension of $\dim_k \Hom( F, E)$ for rank
    $r$, determinant $L$ vector bundles $E$ on $C$ is a multiple of $\hcf( r, d)$.
  \end{itemize}
  (ii follows from \ref{rank_x_weight} and example \ref{Hombdl}.ii. For $F$ generic, it also follows from
  the theorem of Hirschowitz \cite{hirschowitz} that $\Hom(F,E) = 0$ or $\Ext^1(F,E) = 0$ for $E,F$ generic.)
\end{example}

\begin{example}
  Given a closed subscheme $\emptyset \neq D \varsubsetneq C$, we consider rank $r$ vector bundles $E$ on $C$
  with $\det( E) \cong L$, endowed with a \emph{level structure} at $D$, i.\,e.\ with a nonzero morphism $E
  \to \O_D^r$. Let $\Levcoarse_{r, L}^D$ be the coarse moduli scheme of such $E$ with level structure
  which are stable in the sense of \cite{seshadri}. Then $\Levcoarse_{r, L}^D$ is rational.
\end{example}
\begin{proof}
  Let $\stackV$ be the vector bundle of weight $-1$ on $\Bunstack_{r, L}$ whose fibre at any point $[E]$ is
  $\Hom_{\O_C}(E, \O_D^r)$. Then $\Levcoarse_{r, L}^D$ is the coarse space associated to an open substack
  of $\Gr_1( \stackV)$; hence it is rational by corollary \ref{rationality}.ii above.
\end{proof}
More generally, the above results can be applied in a similar way to \emph{stable pairs} on $C$ in the sense
of Huybrechts and Lehn \cite{huybrechts-lehn_pairs}, and to other sorts of \emph{augmented bundles} on $C$ in
the sense of \cite{BDGW}; we mention:
\begin{example}
  A \emph{coherent system} $(E, V)$ of type $(r, d, n)$ consists of a rank $r$, degree $d$ vector bundle $E$
  on $C$ and an $n$-dimensional vector subspace $V \subseteq \cohom^0( E)$, cf. \cite{lepotier},
  \cite{raghavendra-vishwanath} or \cite{king-newstead}. Assuming $k = \complexnums$, $n < r$, $0 < d$ and
  $n \leq d+(r-n)(g-1)$, let $\Cohcoarse_{r, d}^n$ be the coarse moduli scheme of such coherent systems which
  are $\alpha$-stable for the largest possible value $\alpha = d/(r-n) -\epsilon$ of the stability parameter
  $\alpha \in \reals$ \cite[\S 3]{BDGW}. Then $\Cohcoarse_{r, d}^n$ is birational to $\Aa^s \times
  \Buncoarse_{h, 0}$ for $h := \hcf( r, d, n)$ and some $s \in \naturals$.
\end{example}
\begin{proof}
  Let $\stackV$ be the vector bundle of weight $-1$ on $\Bunstack_{r-n, d}^{\stab}$ whose fibre at any
  point $[E']$ is $\Ext^1( E', \O_C)$; note that $\Hom( E', \O_C) = 0$ due to the stability of $E'$, so
  $\stackV$ is indeed a vector bundle of rank $d + (r-n)(g-1) \geq n$. The coarse space of the Grassmannian
  $\Gr_n( \stackV)$ over $\Bunstack_{r-n, d}^{\stab}$ is birational to $\Cohcoarse_{r, d}^n$ according to the
  proof of \cite[Thm. 4.5]{bradlow-garciaprada}; now apply corollary \ref{rationality}.i over the function
  field of $\Pic^d( C)$.
\end{proof}

\begin{example}
  Let $\BNcoarse^n_{r, d} \subseteq \Buncoarse_{r, d}$ be the \emph{Brill-Noether locus} of vector bundles
  $E$ on $C$ with $\dim_k \cohom^0(E) \geq n$. Assume $n=1$ or that $C$ is a Petri curve and $n \in \{2,3\}$.
  If $k = \complexnums$, $n<r$ and $0<d < r(g-1)+n$, then $\BNcoarse^n_{r,d}$ is birational to $\Aa^s \times
  \Buncoarse_{h,0}$ for $h := \hcf(r,d,n)$ and some $s \in \naturals$ (or empty if $n=3$, $r=4$, $d=1$ and
  $g=2$).
\end{example}
\begin{proof}
  Here $\BNcoarse_{r, d}^n$ is birational to $\Cohcoarse_{r, d}^n$, according to \cite[Section 11.2]{BGMN}.
\end{proof}

There seems to be no easy generalisation of theorem \ref{mainthm} to weights $w \notin \{0, \pm 1\}$, mainly
because an analogue of the Hecke correspondence \ref{hecke} is missing; cf. also remark \ref{weight_x}. But
at least, we have the following:
\begin{cor}
  Let $L$ be a line bundle of degree $d$ on $C$, let $\stackV$ be a vector bundle of weight $w$ on an open
  substack $\emptyset \neq \stackU \subseteq \Bunstack_{r, L}$, and let $j \leq \rank( \stackV)$. If $r$, $d$
  and $w \cdot j$ have no common factor $\neq \pm 1$, then $\Gr_j( \stackV)$ is stably birational to
  $\Bunstack_{1, L} \cong \BGm$.
\end{cor}
\begin{proof}
  Replacing $r$ by $\hcf(r, d)$ and $k$ by a rational function field over $k$, theorem \ref{mainthm} allows
  us to assume $d=0$ without loss of generality. This case follows from proposition \ref{weightreduction}.iii
  and theorem \ref{mainthm} again.
\end{proof}

\begin{example}
  Over $k = \complexnums$, let $\rho: \GL( r) \to \GL( V)$ be an algebraic representation that maps the
  center $\Gm \subseteq \GL( r)$ to the center $\Gm \subseteq \GL( V)$ with weight $w \in \integers$. Let
  $M$ be a line bundle on $C$. We consider vector bundles $E$ of rank $r$ and degree $d$ on $C$ which are
  \emph{decorated} in the sense of A. Schmitt \cite{schmitt}, i.\,e.\ endowed with an element of $\mathbb{P
  }\Hom( E_{\rho}, M)$ where $E_{\rho}$ is the vector bundle on $C$ with fibre $V$ associated to the $\GL(r)
  $-bundle given by $E$.

  We fix a stability parameter $\delta > 0$ and choose $M$ of sufficiently large degree such that the
  resulting coarse moduli scheme $\Deccoarse_{r, d}^M$ of $\delta$-stable decorated vector bundles is
  irreducible by \cite[Thm. 3.5]{schmitt}. If $r, d$ and $w$ have no common factor $\neq \pm 1$, then the
  generic fibre of $\det: \Deccoarse_{r, d}^M \to \Pic^d( C)$ is stably rational.
\end{example}


\begin{thebibliography}{10}

\bibitem{amitsur}
S.~Amitsur.
\newblock Generic splitting fields of central simple algebras.
\newblock {\em Ann. of Math. (2)}, 62:8--43, 1955.

\bibitem{artin}
M.~Artin.
\newblock Brauer-{S}everi varieties.
\newblock In {\em Brauer groups in ring theory and algebraic geometry (Antwerp
  1981)}, pages 194--210. Springer-Verlag, Berlin, 1982.

\bibitem{beauville}
A.~Beauville.
\newblock Fibr\'es de rang {$2$} sur une courbe, fibr\'e d\'eterminant et
  fonctions th\^eta.
\newblock {\em Bull. Soc. Math. France}, 116(4):431--448 (1989), 1988.

\bibitem{blanchet}
A.~Blanchet.
\newblock Function fields of generalized {B}rauer-{S}everi varieties.
\newblock {\em Comm. Algebra}, 19(1):97--118, 1991.

\bibitem{boden-yokogawa}
H.~U. Boden and K.~Yokogawa.
\newblock Rationality of moduli spaces of parabolic bundles.
\newblock {\em J. London Math. Soc. (2)}, 59(2):461--478, 1999.

\bibitem{bourbaki_AC}
N.~Bourbaki.
\newblock {\em Commutative algebra. {C}hapters 1--7}.
\newblock Elements of Mathematics. Springer-Verlag, Berlin, 1989.

\bibitem{BDGW}
S.~Bradlow, G.~Daskalopoulos, O.~Garc{\'{\i}}a-Prada, and R.~Wentworth.
\newblock Stable augmented bundles over {R}iemann surfaces.
\newblock In {\em Vector bundles in algebraic geometry (Durham 1993)}, pages
  15--67. Cambridge Univ. Press, Cambridge, 1995.

\bibitem{bradlow-garciaprada}
S.~Bradlow and O.~Garc{\'{\i}}a-Prada.
\newblock An application of coherent systems to a {B}rill-{N}oether problem.
\newblock {\em J. Reine Angew. Math.}, 551:123--143, 2002.

\bibitem{BGMN}
S.~Bradlow, O.~Garc{\'{\i}}a-Prada, V.~Mu{\~n}oz, and P.~Newstead.
\newblock Coherent systems and {B}rill-{N}oether theory.
\newblock {\em Internat. J. Math.}, 14(7):683--733, 2003.

\bibitem{faddeev}
D.~K. Faddeev.
\newblock The group of divisor classes on some algebraic curves.
\newblock {\em Soviet Math. Dokl.}, 2:67--69, 1961.

\bibitem{faltings}
G.~Faltings.
\newblock Vector bundles on curves.
\newblock Lectures held in Bonn, 1995.

\bibitem{giraud}
J.~Giraud.
\newblock {\em Cohomologie non ab\'elienne}.
\newblock Springer-Verlag, Berlin, 1971.

\bibitem{ega3}
A.~Grothendieck.
\newblock {EGA III}: \'{E}tude cohomologique des faisceaux coh\'erents.
\newblock {\em Inst. Hautes \'Etudes Sci. Publ. Math.}, 11/17, 1961/63.

\bibitem{ega4}
A.~Grothendieck.
\newblock {EGA IV}: \'{E}tude locale des sch\'emas et des morphismes de
  sch\'emas.
\newblock {\em Inst. Hautes \'Etudes Sci. Publ. Math.}, 20/24/28/32, 1964--67.

\bibitem{sga4}
A.~Grothendieck et~al.
\newblock {\em {SGA IV}: Th\'eorie des topos et cohomologie \'etale des
  sch\'emas}, volume 269, 270, 305 of {\em Lecture Notes in Mathematics}.
\newblock Springer-Verlag, Berlin, 1972/73.

\bibitem{jochen}
J.~Heinloth.
\newblock {\"U}ber den {M}odulstack der {V}ektorb\"undel auf {K}urven.
\newblock Diploma thesis, University of Bonn, 1998.

\bibitem{hirschowitz}
A.~Hirschowitz.
\newblock {P}robl\`emes de {B}rill-{N}oether en rang sup\`erieur.
\newblock \\ \verb|http://math.unice.fr/~ah/math/Brill/|.

\bibitem{rat}
N.~Hoffmann.
\newblock Moduli stacks of vector bundles on curves and the {K}ing-{S}chofield
  rationality proof.
\newblock preprint, available at www.arXiv.org.

\bibitem{huybrechts-lehn_pairs}
D.~Huybrechts and M.~Lehn.
\newblock Stable pairs on curves and surfaces.
\newblock {\em J. Algebraic Geom.}, 4(1):67--104, 1995.

\bibitem{huybrechts-lehn}
D.~Huybrechts and M.~Lehn.
\newblock {\em The geometry of moduli spaces of sheaves}.
\newblock Friedr. Vieweg \& Sohn, Braunschweig, 1997.

\bibitem{king-schofield}
A.~King and A.~Schofield.
\newblock Rationality of moduli of vector bundles on curves.
\newblock {\em Indag. Math. (N.S.)}, 10(4):519--535, 1999.

\bibitem{king-newstead}
A.~D. King and P.~E. Newstead.
\newblock Moduli of {B}rill-{N}oether pairs on algebraic curves.
\newblock {\em Internat. J. Math.}, 6(5):733--748, 1995.

\bibitem{Pic_G}
F.~Knop, H.~Kraft, and T.~Vust.
\newblock The {P}icard group of a {$G$}-variety.
\newblock In {\em Algebraische Transformationsgruppen und Invariantentheorie},
  pages 77--87. Birkh\"auser, Basel, 1989.

\bibitem{knudsen-mumford}
F.~F. Knudsen and D.~Mumford.
\newblock The projectivity of the moduli space of stable curves. {I}.
  {P}reliminaries on ``det'' and ``{D}iv''.
\newblock {\em Math. Scand.}, 39(1):19--55, 1976.

\bibitem{knutson}
D.~Knutson.
\newblock {\em Algebraic spaces}.
\newblock Springer-Verlag, Berlin, 1971.
\newblock Lecture Notes in Mathematics, Vol. 203.

\bibitem{laumon}
G.~Laumon and L.~Moret-Bailly.
\newblock {\em Champs alg\'ebriques}.
\newblock Springer-Verlag, Berlin, 2000.

\bibitem{lepotier}
J.~Le~Potier.
\newblock {\em Syst\`emes coh\'erents et structures de niveau}, volume 214 of
  {\em Ast\'erisque}.
\newblock Soci\'et\'e Math\'ematique de France, Paris, 1993.

\bibitem{lieblich}
M.~Lieblich.
\newblock {Moduli of twisted sheaves}.
\newblock preprint math.AG/0411337.

\bibitem{lieblich2}
M.~Lieblich.
\newblock {Twisted sheaves and the period-index problem}.
\newblock preprint math.AG/0511244.

\bibitem{mehta-seshadri}
V.~B. Mehta and C.~S. Seshadri.
\newblock Moduli of vector bundles on curves with parabolic structures.
\newblock {\em Math. Ann.}, 248(3):205--239, 1980.

\bibitem{narasimhan-ramanan}
M.~S. Narasimhan and S.~Ramanan.
\newblock Moduli of vector bundles on a compact {R}iemann surface.
\newblock {\em Ann. of Math. (2)}, 89:14--51, 1969.

\bibitem{newstead}
P.~E. Newstead.
\newblock {\em Introduction to moduli problems and orbit spaces}, volume~51 of
  {\em TIFR Lectures on Mathematics and Physics}.
\newblock Tata Institute of Fundamental Research, Bombay, 1978.

\bibitem{raghavendra-vishwanath}
N.~Raghavendra and P.~Vishwanath.
\newblock Moduli of pairs and generalized theta divisors.
\newblock {\em Tohoku Math. J. (2)}, 46(3):321--340, 1994.

\bibitem{ramanan}
S.~Ramanan.
\newblock The moduli spaces of vector bundles over an algebraic curve.
\newblock {\em Math. Ann.}, 200:69--84, 1973.

\bibitem{russo-teixidor}
B.~Russo and M.~Teixidor~i Bigas.
\newblock On a conjecture of {L}ange.
\newblock {\em J. Algebraic Geom.}, 8(3):483--496, 1999.

\bibitem{schmitt}
A.~Schmitt.
\newblock {A universal construction for moduli spaces of decorated vector
  bundles over curves}.
\newblock {\em Transform. Groups}, 9:167--209, 2004.

\bibitem{schofield-vandenbergh}
A.~Schofield and M.~Van~den Bergh.
\newblock The index of a {B}rauer class on a {B}rauer-{S}everi variety.
\newblock {\em Trans. Amer. Math. Soc.}, 333(2):729--739, 1992.

\bibitem{seshadri}
C.~S. Seshadri.
\newblock {\em Fibr\'es vectoriels sur les courbes alg\'ebriques}, volume~96 of
  {\em Ast\'erisque}.
\newblock Soci\'et\'e Math\'ematique de France, Paris, 1982.

\bibitem{wadsworth}
A.~Wadsworth.
\newblock The index reduction formula for generic partial splitting varieties.
\newblock {\em Comm. Algebra}, 21(4):1063--1070, 1993.

\end{thebibliography}

\end{document}